\documentclass[3p,times]{elsarticle}

\usepackage{amsmath}
\usepackage{color}
\usepackage{bm}
\usepackage[linesnumbered]{algorithm2e}
\usepackage{hyperref}

\usepackage{amssymb}
\usepackage{caption}
\usepackage{subcaption}
\DeclareMathOperator\erf{erf}

\biboptions{compress}

\begin{document}

\begin{frontmatter}




\title{A new lattice Boltzmann scheme for linear elastic solids: periodic problems}


\author[1]{Oliver Boolakee}
\author[2]{Martin Geier}
\author[1]{Laura De Lorenzis\corref{cor1}}
\ead{ldelorenzis@ethz.ch}

\address[1]{Department of Mechanical and Process Engineering, ETH Zürich, 8092 Zürich, Switzerland}
\address[2]{Institute for Computational Modeling in Civil Engineering, TU Braunschweig, 38106 Braunschweig, Germany}
\cortext[cor1]{Corresponding author}

\begin{abstract}
We propose a new second-order accurate lattice Boltzmann scheme that solves the quasi-static equations of linear elasticity in two dimensions. In contrast to previous works, our formulation solves for a single distribution function with a standard velocity set and avoids any recourse to finite difference approximations. As a result, all computational benefits of the lattice Boltzmann method can be used to full capacity. The novel scheme is systematically derived using the asymptotic expansion technique and a detailed analysis of the leading-order error behavior is provided. As demonstrated by a linear stability analysis, the method is stable for a very large range of Poisson's ratios. We consider periodic problems to focus on the governing equations and rule out the influence of boundary conditions. The analytical derivations are verified by numerical experiments and convergence studies. 
\end{abstract}


\end{frontmatter}

% \linenumbers

\section{Introduction}

The lattice Boltzmann method \citep{McNamara1988,Chen1998,Succi2001,Kruger2017} is a numerical method that -- in its native variant -- is primarily used to solve fluid dynamics problems for nearly incompressible flows \citep{Lallemand2021}. In contrast to alternative approaches such as the finite difference or the finite element method, the numerical scheme is not obtained by applying some local or non-local approximation to the exact differentials. Rather, the lattice Boltzmann method is motivated by a simplified microscopic model of a gas that is still capable of recovering macroscopic fluid behavior. This gas kinetic legacy comes with some advantageous properties, such as algorithmic simplicity and good scaling during parallelization. For these reasons, the lattice Boltzmann method is well suited to be run efficiently on modern massively parallel computing architectures \citep{Pasquali2016}.

Within a more general interpretation, the method can be viewed from a purely numerical standpoint without any connection to the kinetic theory of gases. This opens up the possibility to consider the lattice Boltzmann method with its already mentioned benefits as a numerical algorithm that can be utilized to find approximate solutions to certain classes of partial differential equations. In this spirit, lattice Boltzmann schemes solving the diffusion equation \citep{Wolf-Gladrow1995}, shallow water equations \citep{Zhou2000}, wave equation \citep{Guangwu2000}, Poisson equation \citep{Chai2008}, conservative phase-field equation \citep{Geier2015} and many more problems have been proposed.

Initial ideas aimed at solving the linear elastodynamic equations with equal longitudinal and transverse wave speeds are reported by Marconi \& Chopard (2003) \citep{Marconi2003}. In this work, a simple lattice Boltzmann scheme solving the wave equation is combined with a velocity-Verlet-type time integration to track the evolution of the displacement field. The distribution function is hereby physically interpreted as the interaction force between neighboring particles and an energy criterion allows to "switch off" bonds that have been loaded beyond a certain threshold. Only qualitative results are shown for a few fracture and fragmentation test cases. The approach by O'Brien et al. (2012) \citep{OBrien2012} solves the wave equations in 2D and 3D using extended velocity sets, but is also limited to equal wave speeds. In order to stabilize the method against oscillations, finite difference schemes are combined with flux limiters known from conventional fluid dynamics methods. A further extension to arbitrary Poisson's ratios and thus different speeds for the transverse and longitudinal waves is suggested by Murthy et al. (2018) \citep{Murthy2018}. However, the constitutive relation assumes hypoelastic material behavior. To solve this problem so-called crystallographic lattices are used in 3D, but the approach still relies on a finite difference approximation to calculate the local volume change. In contrast to the previous works, a mesh study with a quantitative error analysis is provided, which reveals less than first-order convergence under grid refinement. Schlüter et al. (2018) \citep{Schluter2018} propose a method to solve the elastodynamic problem under anti-plane shear deformation, which significantly simplifies the equations. The resulting scalar wave equation is solved using the lattice Boltzmann scheme presented in \citep{Guangwu2000}. Furthermore, their contribution introduces formulations to handle Dirichlet- and Neumann-type boundary conditions so that physically relevant problems can be solved as well. Their method is extended to solve the general equations of elastodynamics \citep{Schluter2021} by applying a split into two wave equations (in 2D) governing the evolution of the dilation and rotation component of the displacement field. In order to combine the results of the two wave equations -- each solved by a lattice Boltzmann scheme -- finite difference approximations for the gradient and rotation operator are used. Along with an extension towards non mesh-conforming boundary formulations, qualitatively reasonable results for the tension test of a plate with hole are obtained.

In contrast to the previous works, Yin et al. (2016) \citep{Yin2016} solve the quasi-static equations of linear elasticity using the lattice Boltzmann method. Because the numerical method can be viewed as an explicit scheme in time, it is unable to solve elliptic problems. For this reason, they propose to extend the target equation by a time-dependent damping term so that the lattice Boltzmann method can iterate on this modified problem until steady state is reached. Following their approach, the solution for the vector-valued displacement field is obtained using multiple distribution functions, with each distribution function governing the evolution of one component of the displacement field. For the evolution of each distribution, the divergence of the displacement field is required as input for the local equilibrium. Although not explicitly stated in the work, it is assumed that this quantity is computed using a finite difference stencil. The contribution introduces two schemes for 2D and 3D problems and convergence studies with simple numerical test cases demonstrate approximately second-order accuracy of the method. Although numerical examples with Dirichlet-type boundary conditions are presented, no information on the boundary formulation is provided.

The main drawbacks of the previously cited studies can be summarized as follows:
\begin{itemize}
    \item All studies use finite difference stencils either to compute spatial gradients or to perform time integration. At first glance this combination of the lattice Boltzmann method with finite difference approximations may not seem problematic, however the high computational efficiency achieved by optimized implementations \citep{Mattila2007,Bailey2009,Geier2017} relies on the very beneficial memory access pattern resulting from the native lattice Boltzmann method. Moreover, the overall order of accuracy is determined by the least accurate step, so that the use of inexpensive low-order finite difference approximations leads to a reduction of the total accuracy. For these reasons, it is generally preferred to avoid computations involving finite difference stencils.
    \item In some studies \citep{OBrien2012,Murthy2018}, the elastodynamic equations are recovered using extended velocity sets, which involve more than only next-neighbor communication. This is expected to reduce the achievable computational efficiency of the implementation. Additionally, when trying to enforce physically consistent boundary conditions, the larger velocity stencils most likely pose significant challenges. Potentially for this reason, the examples shown in \citep{OBrien2012,Murthy2018} involve only periodic domains.
    \item In \citep{Yin2016,Schluter2021}, multiple distribution functions must be solved for to determine the vector-valued solution field. However, for improved numerical efficiency and reduction of memory requirements, solving for a single distribution function is of paramount importance.
\end{itemize}
In summary, so far no lattice Boltzmann scheme can accurately solve the dynamic or quasi-static equations of linear elasticity without recourse to finite difference approximations and/or extended velocity sets or multiple distribution functions. 

In this work, we design a novel second-order consistent lattice Boltzmann scheme solving the quasi-static equations of linear elasticity, which uses a single distribution function with a standard velocity set and does not resort to any finite difference approximations. We demonstrate that the method is stable and accurate for a large range of Poisson's ratios. Within the scope of the present work, we focus on the approximation of the governing equations in the bulk and rule out the influence of boundary conditions by considering periodic problems. All derivations and numerical examples are carried out in 2D in order to keep the expressions manageable.

This manuscript is organized as follows. In the following Section \ref{sec:equation}, the target equation of quasi-static linear elasticity with the extension by a damping term is introduced. The derivation of the novel lattice Boltzmann scheme is outlined in Section \ref{sec:lbm}, followed by an investigation of the leading-order error in Section \ref{sec:conv} and a stability analysis in Section \ref{sec:stability}. Section \ref{sec:icbc} describes how initial conditions and periodic boundary conditions are handled. The final Section \ref{sec:num} shows the results of numerical convergence studies with the purpose of verifying the analytical derivations in the previous sections using the method of manufactured solutions \citep{Roache2002}.

\section{Target problem of linear elasticity} \label{sec:equation}

This section provides a brief introduction to the equations of linear elasticity in 2D. In a next step, the quasi-static equations are endowed with a damping term so that the modified problem can be solved by the lattice Boltzmann method. Finally, we present a consistent non-dimensionalization of the problem.

\subsection{Modified quasi-static linear elasticity in 2D}

The target equation of 2D linear elasticity under the quasi-static assumption is obtained by combining three ingredients: 
\begin{itemize}
\item
the linear momentum balance equation
\begin{equation}
     \nabla \cdot \boldsymbol{\sigma} + \boldsymbol{b} = \boldsymbol{0},
\label{eq:balance}
\end{equation}
where $\boldsymbol{\sigma}$ denotes the second-order Cauchy stress tensor and $\boldsymbol{b}$ the external body load;
\item the linear elastic isotropic constitutive law
\begin{equation}
    \boldsymbol{\sigma}=\lambda \mathrm{tr}(\boldsymbol{\varepsilon}) \boldsymbol{I} + 2\mu \boldsymbol{\varepsilon},
    \label{eq:Hooke}
\end{equation}
where $\boldsymbol{\varepsilon}$ denotes the second-order infinitesimal strain tensor, $\lambda$ and $\mu$ are known as Lam\'{e} parameters ($\mu$ is also called shear modulus), and $\boldsymbol{I}$ is the second-order unit tensor, and
\item the linear kinematic relations
\begin{equation}
    \boldsymbol{\varepsilon} = \mathrm{\frac{1}{2}}\left(\nabla\boldsymbol{u}+ (\nabla\boldsymbol{u})^T\right),
    \label{eq:kin}
\end{equation}
where $\boldsymbol{u}$ is the displacement field. 
\end{itemize}
All three equations are valid in each point of the domain $\Omega\subset \mathbb{R}^2$. Note that we assume this domain to lie in 2D space and not to be a 2D manifold embedded in 3D space, as more common in the mechanics literature, e.\,g. with the plane strain or plane stress assumption. Accordingly, the relations of the Lam\'{e} parameters to Young's modulus $E$ and Poisson's ratio $\nu$ read as follows 
\begin{equation}
    \label{eq:material_relations}
    \lambda = \frac{E \nu}{1-\nu^2}  \qquad\qquad \mu = \frac{E}{2(1+\nu)}\qquad\qquad K = \lambda + \mu = \frac{E}{2(1-\nu)},
\end{equation}
where we also introduced the bulk modulus $K$. The 2D assumption is not critical for the derivation of the numerical scheme, which can be carried out for the plane strain and plane stress cases in an analogous fashion; it is chosen in this paper as it leads to the most straightforward relations.

It is clear from Eq. \eqref{eq:material_relations}$_3$  that in this case incompressible material behavior is obtained for $\nu = 1$. The combination of Eqs. \eqref{eq:balance}-\eqref{eq:kin} delivers the governing equation in the primary variable $\boldsymbol{u}$
\begin{equation}
     \mu \nabla^2 \boldsymbol{u} + K \nabla (\nabla \cdot \boldsymbol{u}) + \boldsymbol{b} = \boldsymbol{0} \qquad \text{in } \Omega \subset \mathbb{R}^2,
\label{eq:linelas}
\end{equation}
also known as the Navier-Cauchy equation.

As proposed by \citep{Yin2016}, the elliptic equation \eqref{eq:linelas} needs to be extended by some time-dependent damping term to be solved by the lattice Boltzmann scheme. To this end, Eq. \eqref{eq:linelas} is replaced by the following time-dependent problem:
\begin{equation}
    \kappa \partial_t \boldsymbol{u} = \mu \nabla^2 \boldsymbol{u} + K \nabla (\nabla \cdot \boldsymbol{u}) + \boldsymbol{b} \qquad \text{in } \Omega  \times [0,t_f] \subset \mathbb{R}^2 \times \mathbb{R}^+_0.
    \label{eq:linelasdamp}
\end{equation}
Differently from \citep{Yin2016}, a damping constant $\kappa$ is also introduced for dimensional consistency and to obtain a more straightforward control over the temporal evolution of the solution. Because of this extension, the domain on which the new relation is defined comprises an open set in space $\Omega$ and a time segment $[0,t_f]$ with final time $t_f$. Assuming that the body load is constant in time, the solution to Eq. \eqref{eq:linelasdamp} evolves towards a steady state as $t_f \rightarrow \infty$. Once steady state is reached,  $\boldsymbol{u}$ naturally fulfills static equilibrium, i.\,e. Eq. \eqref{eq:linelas}. In practice, $t_f$ is chosen large enough such that $\lVert\partial_t \boldsymbol{u}\rVert < \textsc{tol}$ in some norm and for a given tolerance $\textsc{tol}$.

In this initial work, Eq. \eqref{eq:linelasdamp} is solved on a periodic domain. For simplicity, it is assumed that the periodicity is aligned with the primary lattice directions, spanned by the unit vectors $\boldsymbol{e}_x$ and $\boldsymbol{e}_y$. The periodic lengths are assumed to be $L_x$ and $L_y$ so that the solution only needs to be computed on a $L_x \times L_y$ subset of $\mathbb{R}^2$. Therefore, the solution to Eq. \eqref{eq:linelasdamp} is subject to the additional constraint
\begin{equation}
    \label{eq:periodicity}
    \boldsymbol{u}(\boldsymbol{x}+L_x \boldsymbol{e}_x + L_y \boldsymbol{e}_y,t) = \boldsymbol{u}(\boldsymbol{x},t) \qquad \forall \boldsymbol{x} \in \mathbb{R}^2, \forall t \in [0,t_f] .
\end{equation}
Finally, Eq. \eqref{eq:linelasdamp} is furnished with the following initial condition:
\begin{equation}
\label{eq:init}
    \boldsymbol{u} = \boldsymbol{u}_0 \quad\text{in } \Omega \times \{0\}.
\end{equation}

\subsection{Non-dimensionalization}\label{sec:nondim}
%The basic variables in the lattice Boltzmann method are the populations $f_{ij}$ indicating the probability of a particle to move a distance $(i\mathbf{e}_x j\mathbf{e}_y)\Delta x$ in the time $\Delta t$. The actual value of $\Delta x$ and $\Delta t$ is required for scaling the result of the calculation to real word distances and times but it is not required in the execution of the algorithm. Where quantities depending on the physical units of space and time enter the method, using metric units would introduce unnecessary conversions in each time step. It is therefore customary to express all physical quantities in so-called lattice units where $\Delta x=1$ and $\Delta t=1$ and restrict the conversion to and from metric units to the pre- and post-processing. 

For the interpretation of the simulation results independently from the physical units, it is customary to state the governing equation in dimensionless form by scaling with the characteristic length $L$, time $T$ and mass $M$, leading to the respective dimensionless quantities $\tilde{\left(\bullet\right)}$ as follows 

%For an efficient implementation of the lattice Boltzmann method it is crucial to transform the physical problem into so-called lattice units (to be introduced later). To this end, all physical quantities need to be consistently nondimensionalized using a set of reference quantities \cite{Kruger2017}. To nondimensionalize the equations of linear elasticity, all quantities are scaled with a combination of the characteristic length $L$, time $T$ and mass $M$, leading to the respective dimensionless quantities $\tilde{\left(\bullet\right)}$ as follows 
\begin{equation}
\begin{alignedat}{6}
    \boldsymbol{x} &= L \tilde{\boldsymbol{x}} \qquad t = T \tilde{t} \qquad \boldsymbol{u} = U \tilde{\boldsymbol{u}} &&\qquad \Omega &&= L^2\tilde{\Omega} \qquad t_f = T \tilde{t}_f \\[5pt]
    \dfrac{\mu}{\kappa} &= \dfrac{T^{-2}M}{L^{-2}T^{-1}M}\tilde{\mu} = L^2 T^{-1} \tilde{\mu} &&\qquad \dfrac{K}{\kappa} &&= L^2T^{-1} \tilde{K} \\[10pt]
    \dfrac{\boldsymbol{\sigma}}{\kappa} &= \dfrac{L^{-1}T^{-2}MU}{L^{-2}T^{-1}M}\tilde{\boldsymbol{\sigma}} = LT^{-1}U\tilde{\boldsymbol{\sigma}} &&\qquad \dfrac{\boldsymbol{b}}{\kappa} &&= \dfrac{L^{-1}T^{-2}M}{L^{-2}T^{-1}M} \tilde{\boldsymbol{b}} = L T^{-1} \tilde{\boldsymbol{b}}.
\end{alignedat}
\end{equation}

Note that as a result of the division by the damping coefficient $\kappa$ the explicit appearance of the reference mass $M$ is removed. Note also that the displacement field is normalized with the length scale $U$, and because the problem is linear this scaling factor can be chosen arbitrarily.
Introducing the above non-dimensionalization into the governing equation, the periodicity constraint and the initial conditions, i.\,e. Eqs. \eqref{eq:linelasdamp}-\eqref{eq:init}, the dimensionless problem reads
\begin{alignat}{4}
    \label{eq:linnondim}
     \partial_{\tilde{t}} &\tilde{\boldsymbol{u}} = \tilde{\mu} \tilde{\nabla}^2 \tilde{\boldsymbol{u}} + \tilde{K} \tilde{\nabla} (\tilde{\nabla} \cdot \tilde{\boldsymbol{u}}) + L U^{-1} \tilde{\boldsymbol{b}} \quad &&\text{in } \tilde{\Omega}  \times [0,\tilde{t}_f] \\
     \label{eq:initnondim}
     &\tilde{\boldsymbol{u}} = \tilde{\boldsymbol{u}}_0 \qquad &&\text{in } \tilde{\Omega}  \times \{0\} \\
     \text{with } &\tilde{\boldsymbol{u}}(\tilde{\boldsymbol{x}}+\tilde{L}_x \boldsymbol{e}_x + \tilde{L}_y \boldsymbol{e}_y,\tilde{t}) = \tilde{\boldsymbol{u}}(\tilde{\boldsymbol{x}},\tilde{t}) \qquad &&\forall \tilde{\boldsymbol{x}} \in \mathbb{R}^2,
     \label{eq:periodicitynondim}
\end{alignat}
where $\tilde{L}_{x}=L_{x}/L$ and $\tilde{L}_{y}=L_{y}/L$ and the differential operators are dimensionless as well. The dimensionless Cauchy stress is then easily obtained as
\begin{equation}
    \tilde{\boldsymbol{\sigma}}= \tilde{K}(\tilde{\nabla} \cdot \tilde{\boldsymbol{u}}) \boldsymbol{I} +
    \tilde{\mu} \left(\tilde{\nabla}\tilde{\boldsymbol{u}}+ (\tilde{\nabla}\tilde{\boldsymbol{u}})^T-(\tilde{\nabla} \cdot \tilde{\boldsymbol{u}}) \boldsymbol{I}\right).
    \label{eq:stress}
\end{equation}

%Finally, in lattice Boltzmann it is customary to transform all quantities into the lattice unit system. In lattice units, the spacing of the regular square grid and the time step size is one and all lattice nodes at each time step can be addressed by some integer index combination.

\section{Lattice Boltzmann scheme for linear elasticity} \label{sec:lbm}

In this section, the lattice Boltzmann scheme used to solve the target equation is derived. For this purpose, some notation conventions and basic definitions are introduced. %Because the proposed scheme is constructed based on a generalized multiple relaxation time (MRT) collision operator, the introduction directly covers this variant of the lattice Boltzmann method. 
After defining the structure of the numerical scheme, the asymptotic analysis technique is applied to understand how the target equation of linear elasticity can be solved by the scheme. In particular, we carry out a consistent derivation of the relaxation rates and equilibrium moments, which determine the governing equation being solved by the method. Finally, it is shown how the approximate Cauchy stress field can be retrieved from the numerical algorithm with minimal post-processing effort.

\subsection{Basic definitions and notation conventions}

\begin{figure}[htbp]
    \centering
    \includegraphics{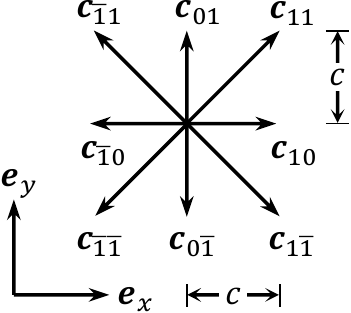}
    \caption{D2Q9 velocity set with the direction of each microscopic velocity $\boldsymbol{c}_{ij}$ indicated by Miller indices, i.\,e. $\bar{1} = -1$}
    \label{fig:lattice}
\end{figure}

The primary solution quantity of the lattice Boltzmann method is the distribution function $f$, thus indicating the connection with the Boltzmann equation (see \citep{Lallemand2021} for an exhaustive review). %Up to an arbitrary multiplicative constant (Laura: I would leave this only in the footnote, here we are not saying what f represents) 
$f$ is a dimensionless function of space $\boldsymbol{x}$, time $t$ and microscopic velocity $\boldsymbol{\xi}$, i.\,e. $f=f(\boldsymbol{x},t,\boldsymbol{\xi})$ \footnote{In the Boltzmann equation, $f$ represents the probability of finding a particle at position $\boldsymbol{x}$ featuring the velocity $\boldsymbol{\xi}$ at a given instant of time $t$. In practice, normalization of the probability distribution to unity is usually omitted such that, strictly speaking, $f$ is a probability density times an arbitrary constant.  }. In the discretized setting of the numerical method however, the distribution function is retained only at a finite number of points in velocity space $\boldsymbol{c}_{ij}$ and the distribution function evaluated at the discrete velocity $\boldsymbol{c}_{ij}$ times some corresponding weight $W_{ij}$ is referred to as the \emph{population} $f_{ij}$, i.\,e.
\begin{equation}
\label{eq:distrib}
    f_{ij}(\boldsymbol{x},t) = W_{ij} f(\boldsymbol{x},t,\boldsymbol{c}_{ij}).
\end{equation}
Moreover, the domain of interest in 2D is discretized by a square lattice with a uniform grid spacing $\Delta x$, whereas the time interval of interest is discretized with a uniform time step $\Delta t$. Adopting the notation of \citep{Geier2015cum}, the two indices $i$ and $j$ denote the unit velocity components in x-direction and in y-direction, respectively; in other words, with the introduction of the scalar lattice speed $c=\Delta x/\Delta t$, the microscopic velocity of the population $f_{ij}$ is defined as $\boldsymbol{c}_{ij}=ic\boldsymbol{e}_x+jc\boldsymbol{e}_y$. Figure \ref{fig:lattice} illustrates the most widely used set of microscopic velocities in 2D, which is known as D2Q9 \citep{Qian1993} as it includes nine velocities on a 2D lattice, i.\,e. $i,j\in\{-1,0,1\}$. Note that the velocity of the so-called rest population $\boldsymbol{c}_{00} = \boldsymbol{0}$ is not visible in this depiction. This is also the set of velocities which will be used in this paper. Hence, the lattice Boltzmann method computes these nine populations $f_{ij}, i,j\in\{-1,0,1\}$ at each point of the lattice for each time step.

\begin{figure}[htbp]
    \centering
    \includegraphics[width=0.8\textwidth]{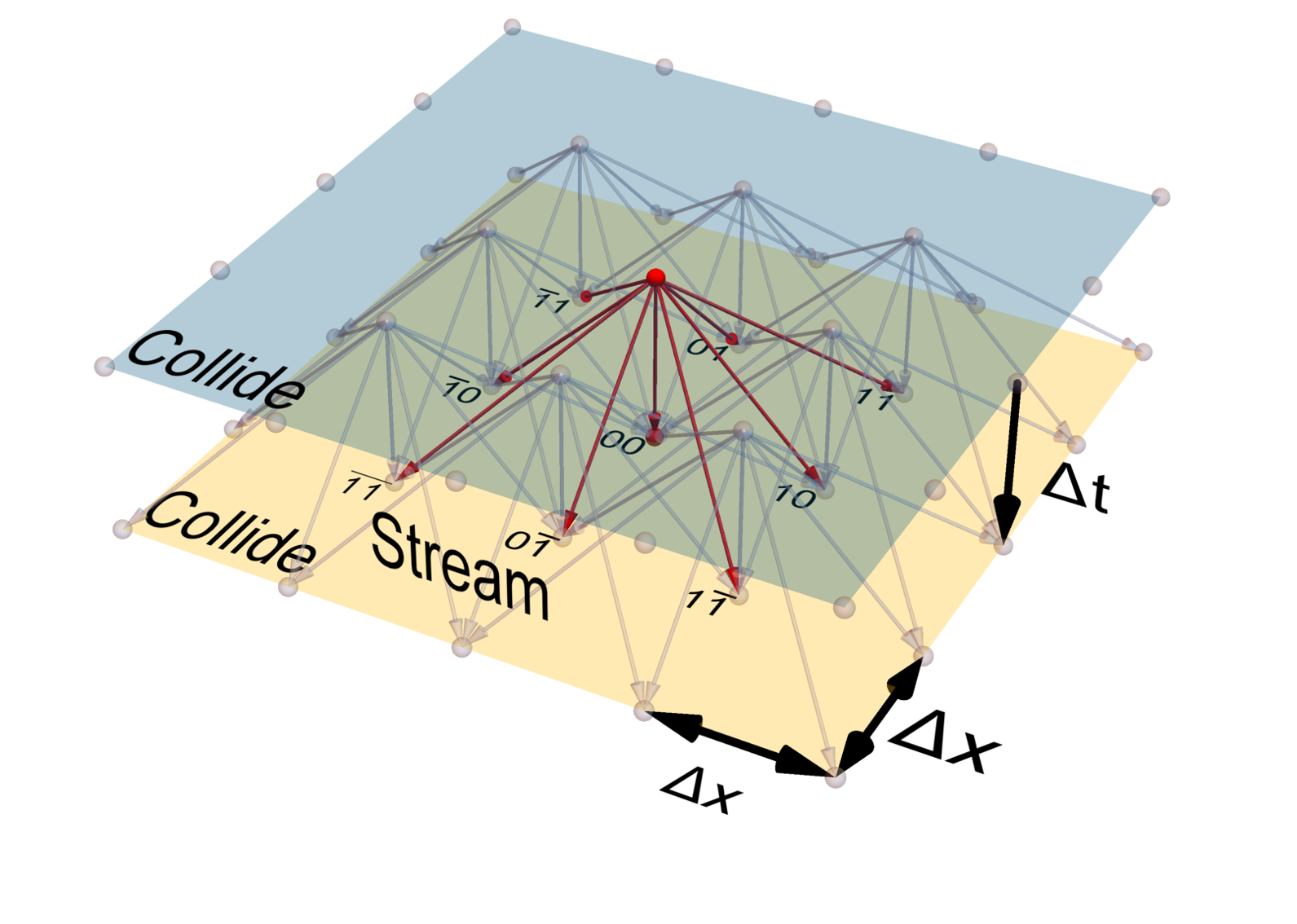}
    \caption{Space-time diagram for the D2Q9 velocity set showing the collision and streaming stages of the central node in red. Between the collision steps the populations move according to their velocity $\boldsymbol{c}_{ij}$. (For  interpretation  of  the  references  to  color  in  this  figure  legend,  the reader  is  referred  to  the  web  version  of  this  article.)}
    \label{fig:stream_collide}
\end{figure}

The computation of the temporal evolution of the populations involves the following two stages in alternation: the collision stage, where the populations $f_{ij}$ locally interact with each other (resulting in a change of their values), and the streaming stage, where each population $f_{ij}$ travels in the direction determined by its microscopic velocity $\boldsymbol{c}_{ij}$ (see Figure \ref{fig:stream_collide}). Both stages will be constructed and described in the following sections. Note that the definition of the discrete set of microscopic velocities on the lattice permits an exact shift of the populations in the streaming stage from one point of the lattice to another without the need for interpolation.

For a given time, once the populations are known at each point of the lattice, it is possible to compute their so-called \emph{raw countable moments} as follows \footnote{These are the discrete and dimensionless counterparts of the statistical moments of the distribution function $f$, defined as integrals over the velocity space.}
\begin{equation}
    m_{ab}=\sum_{i=-1}^1 \sum_{j=-1}^1 i^a j^b f_{ij}.
    \label{eq:moments}
\end{equation}
These are by construction dimensionless. Later on, we will identify some of these moments with the dimensionless solution fields appearing in the governing Eq. \eqref{eq:linnondim} \citep{Geier2015cum,Geier2018}. The so-called order of the moment $m_{ab}$ is given by $a+b$.

It is important to note that the size of the velocity set and thus the number of populations $f_{ij}$ governs the number of independent moments that can be obtained. For example, in the case of the D2Q9 velocity set where $i,j \in \{-1,0,1\}$, it can be easily verified that nine independent raw countable moments exist: $\{m_{00},m_{10},m_{01},m_{11},m_{20},m_{02},m_{12},m_{21},m_{22}\}$, as the computation of all other higher-order moments using Eq. \eqref{eq:moments} leads to expressions that are identical to the ones contained in the set shown above, e.\,g. $m_{30}=m_{10}$ or $m_{14}=m_{12}$.
Because these moments are identified with the macroscopic variables of the physical problem being solved by the method, the size of the velocity set limits how complex the target problem can be. On the other hand, increasing the velocity set, as in the case of the extended velocity sets mentioned in the introduction, comes with drastically increased memory requirements and increased computational cost. Consequently, the velocity set should be chosen only as large as necessary for a given application. For the present study, as we will see later, the standard D2Q9 \citep{Qian1993} velocity set provides enough independent moments to solve the target problem, i.\,e. Eqs. \eqref{eq:linnondim}-\eqref{eq:periodicitynondim}. A significant advantage of this standard velocity set is that it involves only next-neighbor node communication, thus enabling optimized implementation strategies \citep{Mattila2007,Bailey2009,Geier2017}. In Section \ref{sec:asympt} it will be shown that the rest population $f_{00}$ can be removed for the present application, which leads to a further reduction in memory requirements and computational cost.

\subsection{Structure of the lattice Boltzmann method} \label{structureLB}
In this section, we describe the structure of the lattice Boltzmann method used in this work, including the two computational stages of collision and streaming. We employ the standard algorithm \citep{He1997,Chen1998} along with the so-called multiple-relaxation-time (MRT) collision operator \citep{DHumieres2002}. The resulting structure is general and can be used to solve different target problems, but some quantities which appear within this structure will be determined in Section \ref{sec:asympt} in such a way that the numerical scheme solves the target equation \eqref{eq:linnondim} with the desired values of the material parameters and of the applied body loads. The possibility to adjust independently both material parameters of linear elasticity is the reason for the choice of the MRT collision operator. %by independently tuning the different relaxation rates, which will be introduced in the following. %The only exception to this are the second-order moments, where a suitable decomposition of the second-order moments is relaxed instead of the moments directly.

\subsubsection{Collision}

%In order to assign different relaxation rates to different moments, 
As introduced earlier, at the collision stage the populations $f_{ij}$ locally interact with each other, hence post-collision populations are computed from the pre-collision ones. As follows, we outline the structure of this computation using a variant of the MRT collision operator \citep{DHumieres2002}. 

As the name indicates, the MRT collision operator involves multiple relaxation times (or equivalently relaxation rates - quantities to be introduced shortly). This allows to model a physical behavior involving multiple independent problem parameters, such as shear and bulk modulus in the case of linear isotropic elasticity. A needed preliminary step is that the populations are first converted into a suitable set of moments. With respect to the standard MRT scheme introduced in \citep{DHumieres2002}, this paper uses a different moment set $\mathcal{C}$ based on the raw countable moments (see Eq. \eqref{eq:moments}) as defined below
\begin{equation}
\label{eq:collision_set}
    \mathcal{C} = \{m_\alpha \,|\,\alpha \in \mathcal{I}=\{00,11,\text{s},\text{d},12,21,22\}\},
\end{equation}
where $\mathcal{I}$ denotes the corresponding index set. In the moment set the elements with numerical indices are the already introduced raw countable moments, whereas the definition of the moments $m_\text{s}$ and $m_\text{d}$ will be provided below. For each of the moments in $\mathcal{C}$ the collision rule reads
\begin{equation}
    m^*_\alpha = \omega_\alpha m^{eq}_\alpha + (1-\omega_\alpha)m_\alpha \qquad \forall \alpha \in \mathcal{I},
    \label{eq:collision}
\end{equation}
where $m^*_\alpha$ denotes the post-collision value of $m_{\alpha}$. This expression implies a relaxation of each moment in $\mathcal{C}$ towards its local equilibrium value $m^{eq}_\alpha$, where the rate of this relaxation is controlled by $\omega_\alpha$.
%Figure \ref{fig:algorithm} illustrates schematically the collision computations, whereas Algorithm \ref{alg:algorithm} provides the corresponding pseudo-code.
As will become clear in Section \ref{sec:asympt}, the governing equation being solved by the method depends on the choice of the relaxation rates $\omega_\alpha$ and of the local equilibrium moments $m^{eq}_\alpha$; in particular, the material parameters of the target problem can be controlled by adjusting the relaxation rates $\omega_\alpha$ appropriately. 
In Section \ref{3.3.2}, we will derive the concrete expressions of $m^{eq}_\alpha$ and $\omega_\alpha$ to obtain the target equation \eqref{eq:linnondim} with the desired values of the material parameters.

Note that the collision is performed at each lattice node for each time step and involves only local information, which is one of the advantageous properties of the method. To fully characterize the collision stage and to explain why the moments $m_{01}, m_{10}, m_{02}$ and $m_{20}$ are missing in $\mathcal{C}$ whereas "new" indices s and d appear, a few more details are discussed in the following.

\paragraph{Spherical and deviatoric moments}
In the set $\mathcal{C}$ the raw countable second-order moments $m_{02}$ and $m_{20}$ are replaced by the moments
%can be interpreted as the components of the dimensionless second-order moment tensor %of the distribution function with respect to some Cartesian basis, as shown below.
%\begin{equation}
%     \left[
%    \begin{matrix}
%    m_{20} & m_{11} \\ m_{11} & m_{02}
%    \end{matrix}
%    \right]  = \left\{ c^{-2} \sum_{i,j} f_{ij}\boldsymbol{c}_{ij} \otimes \boldsymbol{c}_{ij}\right\}
%\end{equation}
%The spherical part of this tensor is described by the newly introduced moment $m_s$, whereas the deviatoric part is fully characterized by the moment $m_{11}$ together with the newly defined moment $m_d$, as shown below
\begin{equation}
    \label{eq:replace_raw}
    \begin{array}{rcl}
        m_\text{s} &=& m_{20}+m_{02} \\
        m_\text{d} &=& m_{20}-m_{02}
    \end{array}
    \quad\Leftrightarrow\quad
    \begin{array}{rcl}
        m_{20} &=& (m_\text{s}+m_\text{d})/2 \\
        m_{02} &=& (m_\text{s}-m_\text{d})/2.
    \end{array}
\end{equation}
In the new moments $m_\text{s}$ and $m_\text{d}$ the indices refer to spherical and deviatoric, respectively. As will become clear later, this decomposition is convenient to ensure that at the level of the target equation different elastic stiffnesses for volumetric and deviatoric deformation states can be obtained.  
After collision, a back-transformation of $m_\text{s}^*$ and $m_\text{d}^*$ into $m^*_{02}$ and $m^*_{20}$ is carried out using the inverse relations shown in Eq. \eqref{eq:replace_raw}\footnote{This forward and backward transformation involving $m_{02}$ and $m_{20}$ as intermediate step is introduced here only for conceptual clarity and especially to allow for the systematic asymptotic expansion analysis in Section \ref{sec:asympt}. For the implementation of the scheme, it is possible and more efficient to directly transform the populations into the moments contained in the set $\mathcal{C}$ for collision and back to the post-collision populations for streaming.}.

\paragraph{Bared and collision moments}
The second remark involves some additional moment definitions, useful for the following developments. Apart from the pre- and post-collision moments, the so-called bared moments $\bar{m}_{\alpha}$ and collision moments $\Omega_{\alpha}$ with $\alpha \in \mathcal{I}$ are defined as follows
\begin{align}
    \label{eq:bared_moment}
    \bar{m}_{\alpha} &= (m_{\alpha}+m^*_{\alpha})/2 \\
    \label{eq:collision_moment}
    \Omega_{\alpha} &= m^*_{\alpha}-m_{\alpha},
\end{align}
with inverse relations
\begin{align}
    \label{eq:bared_moment_inv}
    \ m_{\alpha}=\bar{m}_{\alpha}-\Omega_{\alpha}/2 \\
    \label{eq:collision_moment_inv}
    \ m^*_{\alpha} = \bar{m}_{\alpha}+\Omega_{\alpha}/2.
\end{align}
These two types of moments are mainly introduced for conceptual purposes, especially during the asymptotic analysis outlined in Section \ref{sec:asympt}, and do not explicitly appear during the actual simulations. The bared moments can be interpreted as the intermediate state of the moments during collision and some of them will be later identified with solution quantities of the target problem. On the other hand, the collision moments are mainly employed to enable the recursive definition of the asymptotic analysis (see Section \ref{sec:asympt}).

For the following developments, it is useful to transform the collision rule in Eq. \eqref{eq:collision} into a version that only involves bared and collision moments using Eqs. \eqref{eq:bared_moment_inv}, \eqref{eq:collision_moment_inv}. The result reads 
\begin{equation}
    \label{eq:gencollision}
    \bar{m}_\alpha = m_\alpha^{eq} - \tau_\alpha \Omega _\alpha\qquad \forall \alpha \in \mathcal{I},
\end{equation}
where the relaxation time $\tau_\alpha$ has been defined as
\begin{equation}
    \label{eq:relaxation_time}
    \tau_\alpha = \frac{1}{\omega_\alpha} - \frac{1}{2},
\end{equation}
with inverse relation
\begin{equation}
    \label{eq:relaxation_rate}
    \omega_\alpha = \left( \tau_\alpha + \frac{1}{2} \right)^{-1}.
\end{equation}

Note that all efficient lattice Boltzmann implementations avoid storing both the pre- and post-collision populations so that it is not advisable to compute the bared moments by their definition in Eq. \eqref{eq:bared_moment}. Moreover, collision moments are never computed, so also Eq. \eqref{eq:gencollision} cannot be used to compute the bared moments. Instead, combining the two versions of the collision rule given in Eqs. \eqref{eq:collision} and \eqref{eq:gencollision} yields
\begin{equation}
    \label{eq:compute_bared}
    \bar{m}_\alpha = \frac{\omega_\alpha}{2}m_\alpha^{eq}+\left(1-\frac{\omega_\alpha}{2} \right) m_\alpha \qquad \forall \alpha \in \mathcal{I}.
\end{equation}
This relation enables to compute the bared moments based only on the pre-collision moments and their equilibrium values, which are both readily accessible during computations.
\paragraph{First-order moments}
Lastly, the special role that the first-order moments $m_{10}$ and $m_{01}$ play during collision is discussed. Since, as will be shown in Section \ref{sec:asympt}, these two moments are identified with the primary variables of the governing equation the scheme is solving for (i.\,e. with the two components of the displacement field), they are not contained in $\mathcal{C}$, because they are not relaxed towards some local equilibrium state. However, during the collision step a forcing is applied to these moments in order to include the body load term of the target problem (see Eq. \eqref{eq:linnondim}). 

\begin{figure}[htbp]
    \centering
    \includegraphics{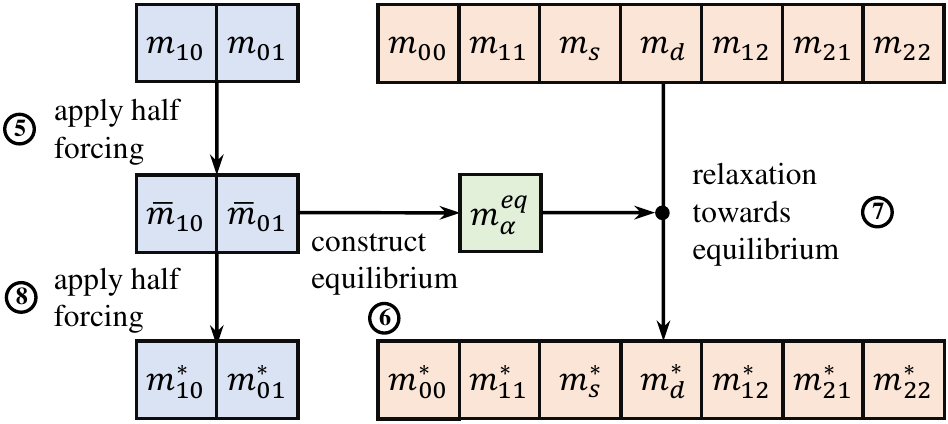}
    \caption{Collision algorithm in moment space (numbers refer to the line number in Algorithm \ref{alg:algorithm})}
    \label{fig:algorithm}
\end{figure}

In accordance with the interpretation of the lattice Boltzmann method as a Strang splitting scheme, the application of the external forcing is split equally in two sub-steps of the collision step \citep{Dellar2013}. The first half of the forcing is applied as follows
\begin{align}
\label{eq:coll10}
    \bar{m}_{10} &= m_{10} + g_x/2 \\
    \label{eq:coll01}
    \bar{m}_{01} &= m_{01} + g_y/2,
\end{align}
to obtain the bared first-order moments, see Figure \ref{fig:algorithm}. Here $g_x$ and $g_y$ represent the x- and y-components of a forcing term that correctly applies the body load appearing in the target equation. The consistent expressions for them will be derived in Section \ref{sec:asympt}. Note that, as will be found in Section \ref{sec:asympt}, some of the local 
equilibrium moments depend on the first-order moments; these local equilibrium moments are computed using the bared first-order moments from \eqref{eq:coll10} and \eqref{eq:coll01}, and not $m_{10}$ and $m_{01}$ (see also Figure \ref{fig:algorithm} and Algorithm \ref{alg:algorithm}). 

After all other moments have undergone the collision as given by Eq. \eqref{eq:collision}, the second half of the forcing term is applied to the bared first-order moments to obtain the post-collision values
\begin{align}
    \label{eq:coll10post}
    m_{10}^* &= \bar{m}_{10} + g_x/2 \\
    \label{eq:coll01post}
    m_{01}^* &= \bar{m}_{01} + g_y/2.
\end{align}

\paragraph{Post-collision populations}
Once all post-collision moments are known, the post-collision populations $f^*_{ij}$, which are required for the next streaming step, are computed by solving 
\begin{equation}
    m^*_{ab}=\sum_{i=-1}^1 \sum_{j=-1}^1 i^a j^b f^*_{ij}
    \label{eq:moments_post}
\end{equation}
for $f^*_{ij}$.

\subsubsection{Streaming}

The second stage of the lattice Boltzmann algorithm is the streaming step, see also Algorithm \ref{alg:algorithm}. It simply involves propagating each post-collision population $f^*_{ij}$, resulting from the collision step, to a neighboring node, where the direction of propagation is given by the microscopic velocity $\boldsymbol{c}_{ij}$, as follows 
\begin{equation}
    \label{eq:streaming}
    f_{ij}(\boldsymbol{x}+\boldsymbol{c}_{ij}\Delta t,t+\Delta t) = f^*_{ij}(\boldsymbol{x},t).
\end{equation}
Because the grid spacing and time step are related to each other by the microscopic lattice speed $c=\Delta x/\Delta t$, streaming implies that the populations "jump" from one lattice node to another, which is one of the reasons for the algorithmic simplicity of the method.

%\subsubsection{Solution algorithm}

%Algorithm \ref{alg:algorithm} summarizes the steps of the numerical solution. 
%The actual value of $\Delta x$ and $\Delta t$ is required for scaling the result of the calculation to real-world distances and times, but it is not required in the execution of the algorithm. Where quantities depending on the physical units of space and time enter the method, using metric units would introduce unnecessary conversions in each time step. %It is therefore customary to express all physical quantities in so-called lattice units where $\Delta x=1$ and $\Delta t=1$ and restrict the conversion to and from metric units to the pre- and post-processing stages. 

\begin{algorithm}[htbp]
\SetKwInOut{Input}{input}\SetKwInOut{Output}{output}
\SetKwBlock{Collide}{Collision}{}
\SetKwBlock{Stream}{Streaming}{}
\Input{populations and forcing at current time step: $f_{ij}(\boldsymbol{x},t)$, $g_x(\boldsymbol{x},t)$ and $g_y(\boldsymbol{x},t)$}
\Output{populations at next time step: $f_{ij}(\boldsymbol{x},t+\Delta t)$, numerical solution of target problem at current time step}
    \For(\tcc*[f]{loop over all nodes}){n $\leftarrow 1$ \KwTo nNodes}{
        \Collide{
            Compute moments \tcc*[f]{Eqs. \eqref{eq:moments} and \eqref{eq:replace_raw}}
            
            Determine solution fields $\tilde{\boldsymbol{u}}(\boldsymbol{x},t)$ and $\tilde{\boldsymbol{\sigma}}(\boldsymbol{x},t)$ \tcc*[f]{see Table \ref{tab:solution}}
            
            Apply half forcing before collision \tcc*[f]{Eqs. \eqref{eq:coll10} and \eqref{eq:coll01}}
            
            Compute local equilibrium moments \tcc*[f]{see Table \ref{tab:scheme}}
            
            Relax all moments in $\mathcal{C}$ \tcc*[f]{Eq. \eqref{eq:collision}}
            
            Apply half forcing after collision \tcc*[f]{Eqs. \eqref{eq:coll10post} and \eqref{eq:coll01post}}
            
            Back-transformation into post-collision populations \tcc*[f]{Eq. \eqref{eq:moments_post}}
        }
        \Stream(\tcc*[f]{Eq. \eqref{eq:streaming}}){
        }
    }
\caption{Single time step of the lattice Boltzmann algorithm}
\label{alg:algorithm}
\end{algorithm}

\subsection{Construction of the scheme using asymptotic analysis}\label{sec:asympt}

The main goal of this section is the systematic determination of the equilibrium moments $m_\alpha^{eq}$, the relaxation rates $\omega_\alpha$, and the forcing components $g_x$ and $g_y$, such that the general lattice Boltzmann scheme in Section \ref{structureLB} (see Figure \ref{fig:algorithm} 
 and Algorithm \ref{alg:algorithm}) solves the target equation \eqref{eq:linnondim} with the desired values of material parameters and body loads. In order to analyze the lattice Boltzmann scheme of Section \ref{structureLB}, we apply the asymptotic expansion technique \citep{Junk2005,Yang2007,Caiazzo2009}, more specifically, the variant presented in \citep{Geier2018}. First, the derivation of the asymptotic expansion from \citep{Geier2018} is briefly outlined. Then we discuss the identification of the equilibrium moments, relaxation rates and forcing terms for the specific target problem of linear elasticity.

\subsubsection{Derivation of the recursive asymptotic expansion}

As follows, we provide an outline of the main steps needed to obtain the recursive asymptotic expansion that will be used in the remainder of this work. The presentation is mainly adopted from \citep{Geier2018}, with slight modifications that do not affect the final expression. The starting point of the asymptotic expansion is the streaming expression in Eq. \eqref{eq:streaming} shifted backwards in time by half a time step:
\begin{equation}
    \label{eq:expansion_step1}
    f_{ij}(\boldsymbol{x}+\boldsymbol{c}_{ij}\Delta t/2,t+\Delta t/2) - f^*_{ij}(\boldsymbol{x}-\boldsymbol{c}_{ij}\Delta t/2,t-\Delta t/2)=0.
\end{equation}
Next, it is assumed that the populations at pre- and post-collision state are regular enough so that we can perform the Taylor series expansion of Eq. \eqref{eq:expansion_step1} around $(\boldsymbol{x},t)$, leading to
\begin{equation}
    \label{eq:expansion_step2}
    \sum_{m,n,l\in\mathbb{N}_0} \frac{1}{m!n!l!} \left(\frac{\Delta t}{2}\right)^m \left(\frac{c \Delta t}{2}\right)^{n+l} \partial_t^m \partial_x^n \partial_y^l \left[ i^n j^l f_{ij} - (-1)^{m+n+l} i^n j^l f_{ij}^* \right] = 0.
\end{equation}
In order to facilitate the physical interpretability of the asymptotic expansion, this expression is transformed into the moment space by pre-multiplying it with $i^a j^b$ and performing the summation over $i$ and $j$. Additionally, the bared and collision moments defined in Eqs. \eqref{eq:bared_moment} and \eqref{eq:collision_moment} are introduced for more concise notation. The result of these algebraic modifications reads
\begin{equation}
    \label{eq:expansion_step3}
    \sum_{m,n,l\in\mathbb{N}_0} \frac{1}{m!n!l!} \left(\frac{\Delta t}{2}\right)^m \left(\frac{c \Delta t}{2}\right)^{n+l} \partial_t^m \partial_x^n \partial_y^l \left[\left(1-(-1)^{m+n+l}\right) \bar{m}_{(a+n)(b+l)} - \frac{1+(-1)^{m+n+l}}{2} \Omega_{(a+n)(b+l)}  \right] = 0.
\end{equation}

Subsequently, the expression is fully nondimensionalized by scaling the time and spatial derivatives with the reference length $L$ and reference time $T$, which were introduced in Section \ref{sec:equation}. In order to obtain only a single smallness parameter $\varepsilon>0$ that needs to be considered during the asymptotic expansion, a relation between $\Delta t$ and $\Delta x=c\Delta t$ is introduced. Because the target equation has a similar structure as the diffusion equation, the uniform grid spacing $\Delta x$ and the time step $\Delta t$ are related according to the so-called diffusive scaling, i.\,e. $\varepsilon^2 \sim (\Delta x)^2 \sim \Delta t$ \citep{Sone1990}. Accordingly, the following expressions are introduced for the grid spacing and time step
\begin{alignat}{6}
    \label{eq:deltax}
    &\Delta x &&= L&&\varepsilon \\
    \label{eq:deltat}
    &\Delta t &&= T&&\varepsilon^2.
\end{alignat}
It will be shown at the end of this section that the diffusive scaling leads to the property that the relaxation rates governing the physical parameters of the target equation do not scale with the smallness parameter. Thus, the limiting process, i.\,e. $\varepsilon \rightarrow 0$, with constant relaxation rates does not alter the physical problem being solved by the method.
Substitution of \eqref{eq:deltax}-\eqref{eq:deltat} into \eqref{eq:expansion_step3} yields the following intermediate result
\begin{equation}
    \label{eq:expansion_step4}
    \sum_{m,n,l\in\mathbb{N}_0} \varepsilon^{2m+n+l} \frac{1}{2^{m+n+l}} \frac{1}{m!n!l!} \partial_{\tilde{t}}^m \partial_{\tilde{x}}^n \partial_{\tilde{y}}^l \left[\left(1-(-1)^{m+n+l}\right) \bar{m}_{(a+n)(b+l)} - \frac{1+(-1)^{m+n+l}}{2} \Omega_{(a+n)(b+l)}  \right] = 0.
\end{equation}

Next, a regular expansion ansatz involving the same smallness parameter $\varepsilon$ is assumed to hold for the populations at pre- and post-collision state, which naturally also applies to all types of moments through their definitions (see Eqs. \eqref{eq:moments}, \eqref{eq:moments_post}, \eqref{eq:bared_moment} and \eqref{eq:collision_moment})
\begin{equation}
    f_{ij} = \sum_{q\in\mathbb{N}_0} \varepsilon^q f_{ij}^{(q)},\quad m_{ab} = \sum_{q\in\mathbb{N}_0} \varepsilon^q m_{ab}^{(q)},\quad m^*_{ab} = \sum_{q\in\mathbb{N}_0} \varepsilon^q m_{ab}^{(q)*},\quad \bar{m}_{ab} = \sum_{q\in\mathbb{N}_0} \varepsilon^q \bar{m}_{ab}^{(q)},\quad \Omega_{ab} = \sum_{q\in\mathbb{N}_0} \varepsilon^q \Omega_{ab}^{(q)}.
    \label{eq:regularExpansion}
\end{equation}
Clearly, all higher-order expansion coefficients vanish during the limiting process $\varepsilon \rightarrow 0$. Therefore, if the method is consistent, the zeroth-order expansion coefficient, i.\,e. that for $q=0$, contains the exact solution and all higher-order terms constitute the numerical error. In order to increase the consistency order to e.\,g. second order, the expansion coefficient with $q=1$ needs to vanish, so that the error scales with $\varepsilon^2$. Thus, the application of the asymptotic expansion derived in this subsection will be used to establish consistency. The following Section \ref{sec:conv} will identify conditions so that the first-order error contributions vanish.

Introducing the regular expansion ansatz into Eq. \eqref{eq:expansion_step4} results in
\begin{equation}
    \label{eq:expansion_step5}
    \sum_{m,n,l,q\in\mathbb{N}_0} \varepsilon^{2m+n+l+q} \frac{1}{2^{m+n+l}}\frac{1}{m!n!l!} \partial_{\tilde{t}}^m \partial_{\tilde{x}}^n \partial_{\tilde{y}}^l \left[\left(1-(-1)^{m+n+l}\right) \bar{m}^{(q)}_{(a+n)(b+l)} - \frac{1+(-1)^{m+n+l}}{2} \Omega^{(q)}_{(a+n)(b+l)}  \right] = 0.
\end{equation}
Because $\varepsilon>0$ is arbitrary, Eq. \eqref{eq:expansion_step5} holds only if the expressions vanish individually at each order in the smallness parameter. Therefore, in order to investigate the relation at some arbitrary order $r$, the summation in Eq. \eqref{eq:expansion_step5} is evaluated with the additional condition $2m+n+l+q=r$. Continuing with the expression obtained by this evaluation, the final expression is established by isolating the $r$th-order expansion coefficient of the collision moment on the other side of the equation by requiring $m+n+l\neq0$ under the sum \citep{Geier2018}:
\begin{equation}
    \Omega_{ab}^{(r)} = \sum_{\substack{m,n,l \in \mathbb{N}_0, \\ m+n+l \neq 0, \\ 2m+n+l+q = r}} \frac{1}{2^{m+n+l}}\frac{1}{m!n!l!} \partial_{\tilde{t}}^m \partial_{\tilde{x}}^n \partial_{\tilde{y}}^l \left[ \left( 1- (-1)^{m+n+l} \right)  \bar{m}_{(a+n)(b+l)}^{(q)} -\left( \frac{1+ (-1)^{m+n+l}}{2} \right)  \Omega_{(a+n)(b+l)}^{(q)} \right].
    \label{eq:expansion}
\end{equation}

Combining Eq. \eqref{eq:expansion} with the collision rule in Eq. \eqref{eq:gencollision} rewritten for each expansion order $q$, we obtain a recursive definition for $\Omega_{ab}^{(r)}$  with arbitrary $a,b,r \in \mathbb{N}_0$ that terminates with $\Omega_{ab}^{(0)}=0$. The final expression obtained by this recursive definition involves only equilibrium moments $m_\alpha^{(q)eq}$ and relaxation times $\tau_\alpha$, which will be chosen in the next subsection such that the target Eq. \eqref{eq:linnondim} is solved by the method in the bulk.

In order to properly carry out the collisions of the second-order moments that are split into their spherical and deviatoric parts, the following additional step is employed: whenever $\bar{m}_{20}^{(q)}$ or $\bar{m}_{02}^{(q)}$ appears with arbitrary $q$, the moment is replaced by its volumetric and deviatoric components using the relations in Eq. \eqref{eq:replace_raw}. Afterwards, the collision rule in Eq. \eqref{eq:gencollision} is applied, followed by a back-transformation into raw countable moments. With this additional step, Eq. \eqref{eq:expansion} can be used for the entire analysis, even if other than raw countable moments are involved in the collision. A few simple examples of Eq. \eqref{eq:expansion} being evaluated are provided in \ref{sec:expansion_examples}.

\subsubsection{Identification of the equilibrium moments, relaxation rates and forcing terms}\label{3.3.2}

An important design step during the construction of a new lattice Boltzmann scheme is the identification of some of the moments with physical quantities of the target equation. In this case the components of the dimensionless displacement field $\tilde{\boldsymbol{u}}$ are identified with the first-order moments. Because there are always as many first-order moments as there are space dimensions considered by the lattice Boltzmann velocity set, all components of the vector-valued displacement field solution can be recovered using a single distribution function, which was one of the major goals outlined in the introduction. Combining this choice with the asymptotic expansion of the lattice Boltzmann scheme leads to
\begin{align}
    \label{eq:ident10}
    \bar{m}_{10}^{(q)} &= u_x^{(q)} \\
    \label{eq:ident01}
    \bar{m}_{01}^{(q)} &= u_y^{(q)},
\end{align}
where $u_x^{(q)}$ and $u_y^{(q)}$ denote the $q$th-order expansion coefficient of the approximate dimensionless displacement field components (for which we do not use the superposed $\,\tilde{}\,$ symbol to avoid overloading the notation). Note that the bared first-order moments, as introduced in Eqs. \eqref{eq:coll10} and \eqref{eq:coll01}, are used for the identification.

In the following, using the asymptotic expansion, we determine the generic form of the governing equations of the first-order moments that results from the lattice Boltzmann method with the MRT collision operator. In a subsequent step,  conditions for the equilibrium moments and relaxation rates are identified so that the target Eq. \eqref{eq:linnondim} is recovered. This is demonstrated here by investigating the behavior of the x-component of the approximate dimensionless displacement field as an example. To this end, the recursive definition of the asymptotic expansion (Eq. \eqref{eq:expansion}) is evaluated with $a=1$ and $b=0$. An analogous derivation can be performed for the y-component as well. However, this does not yield any new information and is therefore not shown here. Evaluating Eq. \eqref{eq:expansion} at the zeroth and first order, i.\,e. $r=0$ and $r=1$, yields the following expressions
\begin{alignat}{4}
    &r=0\text{:} \qquad \Omega_{10}^{(0)}&&=0 \\
    \label{eq:first_order_expansion}
    &r=1\text{:} \qquad \Omega_{10}^{(1)}&&=\frac{1}{2}\partial_{\tilde{x}}m_\text{s}^{(0)eq} + \frac{1}{2}\partial_{\tilde{x}}m_\text{d}^{(0)eq} + \partial_{\tilde{y}}m_{11}^{(0)eq}.
\end{alignat}
These relations do not include any quantity involved in the target equations and the only guidance they provide in designing the method is that they have to be fulfilled. This is the case if $\Omega_{10}^{(0)}=\Omega_{10}^{(1)}=0$, which requires that the equilibrium second-order moments are constant to zeroth order in $\varepsilon$.

The next higher-order expression for the first-order moment yields an expression describing the leading-order equivalent partial differential equation being solved by the generic scheme. To this end, the recursion is evaluated with $r=2$. Furthermore, the left hand-side of Eq. \eqref{eq:expansion} is expanded using Eqs. \eqref{eq:coll10} and \eqref{eq:coll10post}, i.\,e. 
\begin{equation}
\Omega_{10}^{(2)} = m_{10}^{(2)*} - m_{10}^{(2)} = g_x^{(2)}.
\end{equation}
Note that $\Omega_{ab}^{(r)}$ is always under our control as it represents the difference between post- and pre-collision moments.
The resulting expression, obtained using the relationships in \ref{sec:expansion_examples}, reads 
\begin{equation}
\label{eq:xcomp_step1}
\begin{split}
    \partial_{\tilde{t}} u_x^{(0)} = &\frac{1}{2}(\tau_\text{s} + \tau_\text{d} ) \partial_{\tilde{x}}^2u_x^{(0)} + \frac{1}{2}(\tau_\text{s} - \tau_\text{d} ) \partial_{\tilde{x}}^2 m_{12}^{(0)eq} + \tau_{11} \partial_{\tilde{y}}^2 m_{12}^{(0)eq} \\ 
    + &\frac{1}{2}(\tau_\text{s} - \tau_\text{d} ) \partial_{\tilde{x}}\partial_{\tilde{y}} u_y^{(0)} + \frac{1}{2}(\tau_\text{s} + \tau_\text{d} + 2\tau_{11})  \partial_{\tilde{x}} \partial_{\tilde{y}} m_{21}^{(0)eq} \\
    - &\frac{1}{2} \partial_{\tilde{x}} m_\text{s}^{(1)eq} - \frac{1}{2} \partial_{\tilde{x}} m_\text{d}^{(1)eq} - \partial_{\tilde{y}} m_{11}^{(1)eq} + g_x^{(2)},
\end{split}
\end{equation}
where the first-order moments have already been replaced using Eqs. \eqref{eq:ident10} and \eqref{eq:ident01}.

Comparing the governing equation of the leading-order expansion coefficient with the x-component of the target equation, i.\,e.
\begin{equation}
\label{eq:targetx}
    \partial_{\tilde{t}} \tilde{u}_x = (\tilde{\mu}+\tilde{K}) \partial_{\tilde{x}}^2 \tilde{u}_x + \tilde{\mu} \partial_{\tilde{y}}^2 \tilde{u}_x + \tilde{K} \partial_{\tilde{x}}\partial_{\tilde{y}} \tilde{u}_y + LU^{-1} \tilde{b}_x,
\end{equation}
immediately reveals a few requirements that the equilibrium moments need to satisfy:
\begin{itemize}
\item 
Firstly, the third-order equilibrium moments need to be proportional to the displacement solution. This is achieved by setting
\begin{align}
    m_{12}^{(q)eq} &= \theta u_x^{(q)} \\
    m_{21}^{(q)eq} &= \theta u_y^{(q)},
\end{align}
where the proportionality constant $\theta$ has been introduced. As will be discussed later, this constant is typically set to a specific value to obtain isotropic behavior of the lattice\footnote{In fluid dynamics, the parameter $\theta$ corresponds to the dimensionless temperature and specifies the thermal velocity of the particles.}.
\item 
Secondly, the first spatial derivatives of the second-order equilibrium moments $m_\text{s}^{(1)eq}, m_\text{d}^{(1)eq}$ and $m_{11}^{(1)eq}$ need to vanish, because no such derivatives appear in the target equation. The most straightforward solution to achieve this is setting 
\begin{equation}
\label{eq:ms_m11}
m_\text{s}^{eq}=m_\text{d}^{eq}=m_{11}^{eq}=0.
\end{equation} 
Note that this choice also trivially satisfies the relation in Eq. \eqref{eq:first_order_expansion}. Introducing these choices into Eq. \eqref{eq:xcomp_step1} leads to
\begin{equation}
    \partial_{\tilde{t}} u_x^{(0)} = \left[ \frac{1}{2}(1-\theta)\tau_\text{d} + \frac{1}{2}(1+\theta)\tau_\text{s} \right] \partial_{\tilde{x}}^2u_x^{(0)} + \theta \tau_{11} \partial_{\tilde{y}}^2 u_x^{(0)} + \left[ \theta \tau_{11} -\frac{1}{2}(1-\theta)\tau_\text{d} + \frac{1}{2}(1+\theta)\tau_\text{s} \right] \partial_{\tilde{x}}\partial_{\tilde{y}} u_y^{(0)} + g_x^{(2)}.
    \label{eq:ux_order0}
\end{equation}
Comparing this expression with Eq. \eqref{eq:targetx} yields the following result for the relaxation times:
\begin{align}
    \label{eq:match11}
    \theta \tau_{11} &= \tilde{\mu} \\
    \label{eq:matchd}
    \frac{1}{2}(1-\theta)\tau_\text{d} &= \tilde{\mu} \\
    \label{eq:matchs}
    \frac{1}{2}(1+\theta)\tau_\text{s} &= \tilde{K}
\end{align}

\item 
Finally, the forcing term of the lattice Boltzmann scheme has to be set as $g_x^{(2)} = LU^{-1} \tilde{b}_x$. 
\end{itemize}
Performing the same steps starting with Eq. \eqref{eq:expansion} and $r=2,a=0,b=1$ leads to analogous results for the equilibrium moments and relaxation rates and to $g_y^{(2)} = LU^{-1} \tilde{b}_y$. 
Altogether, the previous analysis shows that, with the above choices of equilibrium moments, relaxation times and forcing terms, the leading-order expansion coefficient $u_x^{(0)}$ solves the target equation in the bulk of the domain as stated in Eq. \eqref{eq:linnondim}. 

A few final remarks are pointed out concerning the yet unspecified equilibrium moments and relaxation rates, as well as the parameter $\theta$, as follows:

\begin{itemize}
\item The zeroth-order moment $m_{00}$ does not influence the leading-order physics so that its evolution does not need to be tracked. As a result, the so-called rest population $f_{00}$ can be removed, which reduces the memory requirements of the scheme by $1/9$.
\item The equilibrium value of the fourth-order moment $m_{22}$ is yet unknown, because it has no impact on the leading-order solution. The following Section \ref{sec:conv} will show that a second-order consistent scheme is achieved by choosing $m_{22}^{eq}=0$.
\item A straightforward and robust choice for all relaxation rates that are not involved in Eq. \eqref{eq:xcomp_step1}, and hence do not affect the physics of the solution, is to set $\omega_{12}=\omega_{21}=\omega_{22}=1$ \citep{Geier2015cum}, leading to $\tau_{12}=\tau_{21}=\tau_{22}=1/2$. By Eq. \eqref{eq:collision} this leads to the behavior that the related moments are set to their equilibrium value during each collision.
\item It can be shown \citep{Geier2006} that there exist rotations that are rotational invariants of the velocity set, which map the second order moments $m_{11}$ and $m_\text{d}$ onto each other, i.\,e.
\begin{equation}
\begin{split}
    \exists &\boldsymbol{R} \in \text{SO}(2) \quad \text{s.\,t.} \quad \boldsymbol{c}_{ij} = \boldsymbol{R} \boldsymbol{c}_{kl} \quad \forall i,j,k,l \in \{-1,0,1\} \\
    \text{and} \quad &m_\text{d} = \gamma\sum_{i_R,j_R=-1}^1 i_R j_R f_{ij} \qquad \text{with } \left[\begin{matrix} i_R \\ j_R \end{matrix}\right] = \boldsymbol{R} \left[\begin{matrix} i \\ j \end{matrix}\right], \gamma \in \mathbb{R}.
\end{split}
\end{equation}
Therefore, these two moments represent exchangeable physical quantities and should hence, by physical intuition, relax with the same rate. 
As a result, $\omega_{11}=\omega_\text{d}$ and equivalently $\tau_{11}=\tau_\text{d}$ need to hold, which requires setting $\theta=1/3$. However, this condition is not strictly necessary, and an equivalent behavior of the target equations to leading order can be obtained by a different choice of $\theta \in(0,1)$. In this case the relaxation rates $\omega_{11}$ and $\omega_\text{d}$ need to take on distinct values as given by Eqs. \eqref{eq:match11} and \eqref{eq:matchd}. For the remainder of the work, the standard value is assumed, i.\,e. $\theta=1/3$.
\end{itemize}

A summary of the equilibrium moments and relaxation times for all moments involved in the collision is provided in Table \ref{tab:scheme}. 

\begin{table}[htbp]
    \centering
    \caption{Summary of equilibrium moments and relaxation times for linear elasticity}
    \begin{tabular}{l|c c}
         Moment  $m_\alpha \in \mathcal{C}$ & $m_\alpha^{eq}$ & $\tau_\alpha$  \\ \hline
         $m_{00}^\dagger$ & $-$ & $-$ \\
 %        10$^\ddagger$ & $-$ & $-$ \\
 %        01$^\ddagger$ & $-$ & $-$ \\
         $m_{11}$ & $0$ & $\theta^{-1}\tilde{\mu}$ \\ 
         $m_\text{s}$ & $0$ & $2(1+\theta)^{-1}\tilde{K}$ \\
         $m_\text{d}$ & $0$ & $2(1-\theta)^{-1}\tilde{\mu}$ \\
         $m_{12}$ & $\theta \bar{m}_{10}$ & $1/2$ \\ 
         $m_{21}$ & $\theta \bar{m}_{01}$ & $1/2$ \\ 
         $m_{22}$ & $0$ & $1/2$ \\ \hline
         \multicolumn{3}{l}{\footnotesize $ ^\dagger$ This moment can be removed}
    \end{tabular}
    \label{tab:scheme}
\end{table}

To conclude this section, let us discuss an important consequence of the choice of diffusive scaling using Eqs. \eqref{eq:match11}-\eqref{eq:matchs}. Since all three equations lead to analogous conclusions, we will consider Eqs. \eqref{eq:match11}.  Let us transform the right-hand side of Eq. \eqref{eq:match11} back into physical units. To this end, reference length $L$ and reference time $T$ are related to lattice spacing $\Delta x$ and time step size $\Delta t$ %(the so-called lattice units) 
through Eqs. \eqref{eq:deltax} and \eqref{eq:deltat} respectively, obtaining
\begin{equation}
    \label{eq:scaling_example}
    \theta \tau_{11} = \tilde{\mu} = L^{-2}T\frac{\mu}{\kappa} = \frac{\varepsilon^2\Delta t}{\varepsilon^2\Delta x^2} \frac{\mu}{\kappa} = \frac{\Delta t}{\Delta x^2}\frac{\mu}{\kappa}.
\end{equation}
Comparing the first and the last terms in Eq. \eqref{eq:scaling_example} shows that the material parameter in physical units $\mu/\kappa$ is unaffected by the smallness parameter $\varepsilon$ for a constant relaxation time  $\tau_{11}$. This implies that diffusive scaling keeps the target physical problem being solved unchanged during the limiting process $\varepsilon \rightarrow 0$ with constant relaxation rates. This property is of paramount importance because it largely simplifies the asymptotic analysis, which would otherwise require an asymptotic expansion of $\tau_\alpha$ such that Eq. (\ref{eq:gencollision}) would no longer be valid in this simple form. 

\subsection{Cauchy stress solution} \label{sec:stress}

Apart from the displacement field there is generally a great interest in the Cauchy stress field. 
In the following, we demonstrate that the stress components are directly retrieved from the bared second-order moments, with no need for further post-processing calculations. Using the asymptotic analysis in Section \ref{sec:asympt}, the expansion of the second-order moments has the following form
\begin{align}
    -\bar{m}_\text{s} &= \varepsilon (1+\theta)\tau_\text{s} \left(\partial_{\tilde{x}} u_x^{(0)} + \partial_{\tilde{y}} u_y^{(0)}\right) + \varepsilon^2 (1+\theta)\tau_\text{s} \left(\partial_{\tilde{x}} u_x^{(1)} + \partial_{\tilde{y}} u_y^{(1)}\right) + \mathcal{O}(\varepsilon^3) \\
    -\bar{m}_d &= \varepsilon (1-\theta)\tau_\text{d} \left(\partial_{\tilde{x}} u_x^{(0)} - \partial_{\tilde{y}} u_y^{(0)}\right) + \varepsilon^2 (1-\theta)\tau_\text{d} \left(\partial_{\tilde{x}} u_x^{(1)} - \partial_{\tilde{y}} u_y^{(1)}\right) + \mathcal{O}(\varepsilon^3) \\
    -\bar{m}_{11} &= \varepsilon \theta\tau_{11} \left(\partial_{\tilde{y}} u_x^{(0)} + \partial_{\tilde{x}} u_y^{(0)}\right) + \varepsilon^2 \theta\tau_{11} \left(\partial_{\tilde{y}} u_x^{(1)} + \partial_{\tilde{x}} u_y^{(1)}\right) + \mathcal{O}(\varepsilon^3). \label{eq:mbar11}
\end{align}
Note that the bared second-order  moments are efficiently calculated with Eq. \eqref{eq:compute_bared}. A comparison of Eq. \eqref{eq:mbar11} with the Cauchy stress relation in Eq. \eqref{eq:stress} together with the condition for the relaxation rate in Eq. \eqref{eq:match11} reveals that
\begin{equation}
    \tilde{\sigma}_{xy}=\tilde{\sigma}_{yx}=-\bar{m}_{11}^{(1)}.
\end{equation}

In Section \ref{sec:conv} it will be shown that the method can be made second-order consistent in the displacement field solution. This means that $u_x^{(1)}$ and $u_y^{(1)}$ vanish and by Eq. \eqref{eq:mbar11} second-order consistency for the numerical approximation of the shear stress is established as well, i.\,e.
\begin{equation}
    \tilde{\sigma}_{xy}=\tilde{\sigma}_{yx}=-\frac{1}{\varepsilon}\bar{m}_{11} + \mathcal{O}(\varepsilon^2).
\end{equation}
Following analogous steps for the other Cauchy stress components yields the results shown below.
\begin{align}
    \tilde{\sigma}_{xx} &=-\frac{1}{\varepsilon}\frac{(\bar{m}_\text{s}+\bar{m}_\text{d})}{2} + \mathcal{O}(\varepsilon^2) \\
    \tilde{\sigma}_{yy} &=-\frac{1}{\varepsilon}\frac{(\bar{m}_\text{s}-\bar{m}_\text{d})}{2} + \mathcal{O}(\varepsilon^2)
\end{align}
Altogether, this section showed that the components of the Cauchy stress field can be obtained from simple algebraic transformations involving quantities already present during the collision stage of the algorithm (see Algorithm \ref{alg:algorithm}). Another important result is that the numerical solution of the Cauchy stress components is second-order consistent if the displacement solution is second-order consistent as well. In comparison, using the finite element method as discretization scheme for the linear elasticity equation (and assuming sufficient regularity of the exact solution), a linear polynomial ansatz for the approximate solution leads to second-order consistency for the displacement field but only first-order consistency for the stress components.

It is now evident that the spherical-deviatoric decomposition of the dimensionless stress tensor reads 

\begin{equation}
\tilde{\boldsymbol{\sigma}}=\tilde{\boldsymbol{\sigma}}_\text{sph}+\tilde{\boldsymbol{\sigma}}_\text{dev},
\end{equation}
with
\begin{equation}
\tilde{\boldsymbol{\sigma}}_\text{sph}=\frac{1}{2}\mathrm{tr}\left(\tilde{\boldsymbol{\sigma}}\right)\boldsymbol{I}=-\frac{1}{2\varepsilon}\bar{m}_\text{s}\boldsymbol{I}\qquad \left\{\tilde{\boldsymbol{\sigma}}_\text{dev}\right\}=-\frac{1}{\varepsilon}\left[\begin{array}{cc}
\bar{m}_\text{d} & \bar{m}_{11}\\
\bar{m}_{11} & -\bar{m}_\text{d}
\end{array}\right],
\end{equation}
which explains the designation of $\bar{m}_\text{s}$ and $\bar{m}_\text{d}$ as spherical and deviatoric moments, respectively. A summary of all dimensionless solution fields obtained from the lattice Boltzmann scheme is provided in Table \ref{tab:solution}.

\begin{table}[htbp]
    \centering
    \caption{Relation between dimensionless solution quantities and moments of the lattice Boltzmann scheme}
    \begin{tabular}{l|l}
         Quantity & Moment expression \\ \hline
         $\tilde{u}_x$ & $\bar{m}_{10}$\\
         $\tilde{u}_y$ & $\bar{m}_{10}$ \\
         $\tilde{\sigma}_{xx}$ & $-\varepsilon^{-1}(\bar{m}_\text{s}+\bar{m}_\text{d})/2$ \\
         $\tilde{\sigma}_{yy}$ & $-\varepsilon^{-1}(\bar{m}_\text{s}-\bar{m}_\text{d})/2$ \\
         $\tilde{\sigma}_{xy}$ & $-\varepsilon^{-1}\bar{m}_{11}$
    \end{tabular}
    \label{tab:solution}
\end{table}

\section{Leading-order error investigation} \label{sec:conv}

The asymptotic expansion technique can not only be used to design the lattice Boltzmann scheme so that it solves the equations of linear elasticity. Additionally, by continuing the expansion to higher orders, conditions to obtain a higher-order consistent method can be identified. This section will present conditions to remove the first-order error contribution so that a second-order accurate method is achieved. As a result the next higher-order contribution of the expansion takes on the role of the leading-order error. Because this error term can generally not be removed, its influence on the solution accuracy is investigated and measures to keep this error small are introduced.

So far we used the asymptotic expansion terms up to $r=2$. Continuing the expansion to the next order, i.\,e. computing Eq. \eqref{eq:expansion} with $r=3$ and $a=1,b=0$ yields
\begin{align}
\label{eq:ux_order1}
\partial_{\tilde{t}} u_x^{(1)} &= \left[ \frac{1}{2}(1-\theta)\tau_\text{d} + \frac{1}{2}(1+\theta)\tau_\text{s} \right] \partial_{\tilde{x}}^2u_x^{(1)} + \theta \tau_{11} \partial_{\tilde{y}}^2 u_x^{(1)} + \left[ \theta \tau_{11} -\frac{1}{2}(1-\theta)\tau_\text{d} + \frac{1}{2}(1+\theta)\tau_\text{s} \right] \partial_{\tilde{x}}\partial_{\tilde{y}} u_y^{(1)} + r_x^{(3)} \\
\label{eq:first_order_forcing}
\text{with } r_x^{(3)} &= B_1(\tau_\alpha)\partial_{\tilde{x}}^3 m_{22}^{(0)eq}+B_2(\tau_\alpha)\partial_{\tilde{x}}\partial_{\tilde{y}}^2 m_{22}^{(0)eq} \qquad \alpha \in \mathcal{I},
\end{align}
where all relations for the equilibrium moments identified in Section \ref{sec:asympt} have already been introduced. Comparing this equation with Eq. \eqref{eq:ux_order0} reveals that the error coefficient $u_x^{(1)}$ also solves the target problem, only with a different body force term. As indicated by Eq. \eqref{eq:first_order_forcing}, this term involves derivatives of the yet unspecified equilibrium moment $m_{22}^{(0)eq}$ and two constants that depend on the relaxation times $\tau_\alpha$. 

In order to obtain a second-order consistent method, the first-order expansion coefficient $u_x^{(1)}$ needs to vanish. This is achieved if the governing Eq. \eqref{eq:ux_order1} admits a null solution, which requires that the body force term is zero and that the problem is furnished with zero initial and boundary conditions. The former requirement can be easily fulfilled. The $r_x^{(3)}$ term in Eq. \eqref{eq:first_order_forcing} is zero if the fourth-order equilibrium moment $m_{22}^{eq}$ is some arbitrary constant at zeroth order in $\varepsilon$. One possible way to achieve this involves setting $m_{22}^{eq}=0$. To satisfy the latter requirement, the initialization needs to be second-order consistent, which will be discussed in Section \ref{sec:icbc}. The periodic boundary conditions also introduced in Section \ref{sec:icbc} do not influence the consistency order of the method in the bulk. As a result, it can be summarized that the linear elasticity problem of Eq. \eqref{eq:ux_order1} on a periodic domain with zero initial condition and no body load admits a zero solution for $u_x^{(1)}$. Performing the same analysis for the y-component of the displacement solution leads to identical conclusions.

Assuming that the first-order expansion coefficients have been successfully set to zero, the next higher-order coefficients $u_x^{(2)}$ and $u_y^{(2)}$ constitute the leading-order error. The structure and properties of the governing equation for $u_x^{(2)}$ are obtained by evaluating Eq. \eqref{eq:expansion} with $r=4,a=1,b=0$ as shown below.
\begin{align}
\label{eq:ux_order2}
\partial_{\tilde{t}} u_x^{(2)} &= \left[ \frac{1}{2}(1-\theta)\tau_\text{d} + \frac{1}{2}(1+\theta)\tau_\text{s} \right] \partial_{\tilde{x}}^2u_x^{(2)} + \theta \tau_{11} \partial_{\tilde{y}}^2 u_x^{(2)} + \left[ \theta \tau_{11} -\frac{1}{2}(1-\theta)\tau_\text{d} + \frac{1}{2}(1+\theta)\tau_\text{s} \right] \partial_{\tilde{x}}\partial_{\tilde{y}} u_y^{(2)} + r_x^{(4)} \\
\text{with } r_x^{(4)} &= C_1(\tau_\alpha,\theta)\partial_{\tilde{x}}^4 u_x^{(0)}+C_2(\tau_\alpha,\theta) \partial_{\tilde{x}}^3\partial_{\tilde{y}}u_y^{(0)} + C_3(\tau_\alpha,\theta) \partial_{\tilde{x}}^2\partial_{\tilde{y}}^2 u_x^{(0)} + C_4(\tau_\alpha,\theta) \partial_{\tilde{x}}\partial_{\tilde{y}}^3 u_y^{(0)} + C_5(\tau_\alpha,\theta) \partial_{\tilde{y}}^4 u_x^{(0)} \nonumber\\ &+ D_1(\tau_\alpha,\theta) \partial_{\tilde{x}}^2g_x^{(2)} + D_2(\tau_\alpha,\theta) \partial_{\tilde{x}}\partial_{\tilde{y}} g_y^{(2)} + D_3(\tau_\alpha,\theta) \partial_{\tilde{y}}^2 g_x^{(2)} \qquad \alpha \in \mathcal{I}
\label{eq:leadingCoefficients}
\end{align}
The analysis of the equation governing the behavior of $u_y^{(2)}$ leads to analogous results and will be not explicitly shown. Note that in the expression above all time derivatives of $u_x^{(0)}$ and $u_y^{(0)}$ have been replaced using the corresponding governing equations (see Eq. \eqref{eq:ux_order0} for the x-component), which is also the origin of the forcing terms $g_x^{(2)}$ and $g_y^{(2)}$. All coefficients $C_1 \dots C_5$ and $D_1 \dots D_3$ depend on the relaxation times $\tau_\alpha$ and the parameter $\theta$ (see \ref{sec:explicit_constants} for the expressions). Comparing Eq. \eqref{eq:ux_order2} with Eq. \eqref{eq:ux_order0} shows that the governing equation of the second-order error contribution has the same structure as the target problem being solved by the method. This time the body force term $r_x^{(4)}$ is composed of two types of contributions: the fourth-order spatial derivatives of the leading-order solution $u_x^{(0)}$, $u_y^{(0)}$, and the second-order spatial derivatives of the leading-order coefficients of the forcing terms $g_x^{(2)}$, $g_y^{(2)}$, which we know to be related to the physical body load components of the target problem, $b_x$ and $b_y$ (see Section \ref{sec:asympt}).

The latter contribution can be removed. Indeed, $b_x$ and $b_y$ are known and we assume the analytical expressions of their derivatives to be available, so that the terms in $r_x^{(4)}$ and $r_y^{(4)}$ containing $g_x^{(2)}$ and $g_y^{(2)}$ can be compensated for by extending the forcing term with the following higher-order correction
\begin{align}
    \label{eq:force_correction_x}
    g_x &= \varepsilon^2 g_x^{(2)} - \varepsilon^4 \left( D_1(\tau_\alpha,\theta) \partial_{\tilde{x}}^2 g_x^{(2)} + D_2(\tau_\alpha,\theta) \partial_{\tilde{x}}\partial_{\tilde{y}} g_y^{(2)} + D_3(\tau_\alpha,\theta) \partial_{\tilde{y}}^2 g_x^{(2)} \right) \\
    \label{eq:force_correction_y}
    g_y &= \varepsilon^2 g_y^{(2)} - \varepsilon^4 \left( D_3(\tau_\alpha,\theta) \partial_{\tilde{x}}^2 g_y^{(2)} + D_2(\tau_\alpha,\theta) \partial_{\tilde{x}}\partial_{\tilde{y}} g_x^{(2)} + D_1(\tau_\alpha,\theta) \partial_{\tilde{y}}^2 g_y^{(2)} \right).
\end{align}

On the other hand, the displacement field solution can take an arbitrary form and is obviously unknown. Therefore, the first terms in $r_x^{(4)}$ can be removed only if $C_i = 0, i = 1 \dots 5$. As outlined in Section \ref{3.3.2}, out of all the relaxation times $\tau_\alpha, \alpha \in \mathcal{I}$, the ones governing the relaxation of the second-order moments ($\tau_{11}, \tau_{s}$ and $\tau_{d}$) are adjusted to match the dimensionless material parameters of the target problem and the parameter $\theta$ is fixed at its standard value. As will be shown below, the combination of the discretization parameters $\Delta x$, $\Delta t$ and the artificial damping coefficient $\kappa$ leaves one degree of freedom in the choice of the aforementioned relaxation times. Moving on to the relaxation times of the third-order moments, it can be shown that with a similar reasoning as for $\tau_{11}$ and $\tau_d$ (see Section \ref{3.3.2}) the third-order moments cannot be considered to be independent. Therefore, the corresponding relaxation times $\tau_{12}$ and $\tau_{21}$ need to be the same. Summarizing, this leaves in total three independent parameters: one coming from $\tau_\text{s}$, $\tau_\text{d}$ and $\tau_{11}$, another one from $\tau_{12}=\tau_{21}$ and lastly $\tau_{22}$. Accordingly, it cannot be expected that all independent conditions $C_i = 0, i = 1 \dots 5$ are satisfied by some combination of the remaining free parameters, which is necessary to achieve $r_x^{(4)} = 0$ and $r_y^{(4)} = 0$. Thus, there is no way to achieve higher than second-order consistency, and the terms $r_x^{(4)}$ and $r_y^{(4)}$ have a pivotal influence on the behavior of the numerical error in the bulk.

To illustrate this point further, consider the following decomposition of the numerical approximation to the displacement solution obtained by the lattice Boltzmann scheme
\begin{equation}
    \bar{m}_{10} = u_x^{(0)} + \varepsilon^2 u_x^{(2)} + \mathcal{O}(\varepsilon^3).
\end{equation}
Note that an analogous relation is obtained for the y-component as well. In order to reduce the leading-order error contribution $\varepsilon^2 u_x^{(2)}$, there exist two possibilities:
\begin{enumerate}
    \item Decrease $\varepsilon$, i.\,e. perform a mesh and time step refinement obeying the diffusive scaling assumption.
    \item Reduce the magnitude of the solution $u_x^{(2)}$ and $u_y^{(2)}$ to the leading-order error governing equation (see Eq. \eqref{eq:ux_order2} for the x-component).
\end{enumerate}
It is important to keep in mind that reducing the numerical error with the first option comes with a considerable increase in computational effort. Therefore, it is highly advisable to initially exploit any possibility to reduce the error according to the second option. Accordingly, this option is investigated in some more detail in the following.

For periodic problems, the magnitude of the solutions $u_x^{(2)}$ and $u_y^{(2)}$ of the linear governing equations (see Eq. \eqref{eq:ux_order2} for the x-component) is proportional to the magnitude of the body force terms $r_x^{(4)}$ and $r_y^{(4)}$. To this end, $r_x^{(4)}$ and $r_y^{(4)}$ should made as small as possible. This is in turn realized for arbitrary solutions $u_x^{(0)}$ and $u_y^{(0)}$ if $C_i \rightarrow 0, i = 1,\dots,5$ (see Eq. \eqref{eq:leadingCoefficients} for the x-component).

For this initial investigation we decided to keep the search space over which the $C_i$ are minimized fairly manageable. This is achieved by fixing all higher-order relaxation rates $\tau_{12}=\tau_{21}$ and $\tau_{22}$ at their standard values (see Table \ref{tab:scheme}). As a result, only the relaxation rates $\tau_{11}, \tau_\text{s}$ and $\tau_\text{d}$, which are used to adjust the dimensionless parameters $\tilde{\mu}$ and $\tilde{K}$, remain as independent parameters in the expressions for $C_1 \dots C_5$.

\begin{figure}[htbp]
    \centering
    \includegraphics{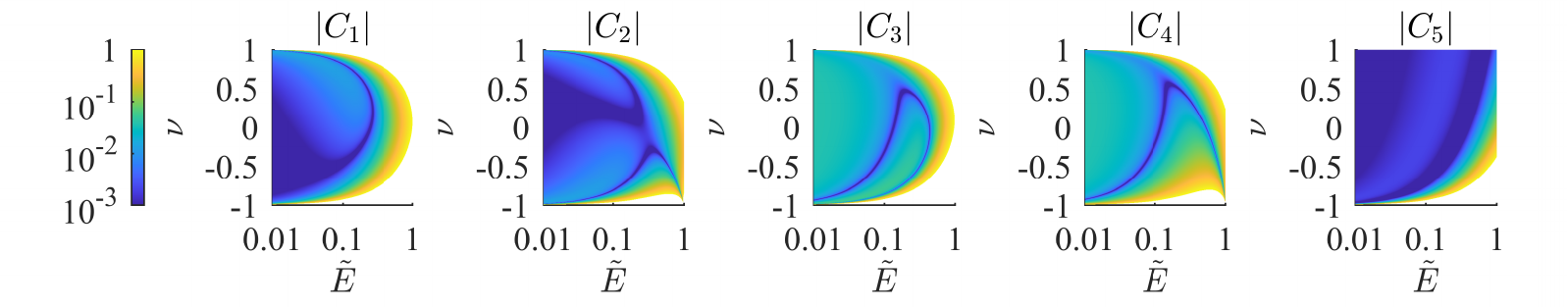}
    \caption[Contour plot of the constants governing the leading-order error]{Contour plot of the constants $C_1 \dots C_5$ governing the leading-order error, see Eq. \eqref{eq:leadingCoefficients}. Values above 1 cut off.}
    \label{fig:errorConstants}
\end{figure}

In order to get a general idea how the body load in the governing equation for the leading-order error can be reduced, $|C_1| \dots |C_5|$ are plotted in Figure \ref{fig:errorConstants} against the dimensionless Young's modulus $\tilde{E}$ and the Poisson's ratio $\nu$ (which are related to $\tilde{\mu}$ and  $\tilde{K}$ by the expressions in Eq. \eqref{eq:material_relations} that also apply to the dimensionless quantities). The dimensionless Young's modulus is related to the physical one by 
\begin{equation}
    \label{eq:scalingE}
    \tilde{E} = L^{-2}T\kappa^{-1}E = \frac{\Delta t}{\kappa \Delta x^2} E,
\end{equation}
see Eqs. \eqref{eq:deltax} and \eqref{eq:deltat}. Eq. \eqref{eq:scalingE} shows that for a given physical Young's modulus $E$ the ratio between time step size $\Delta t$ and grid spacing $\Delta x$ as well as damping constant $\kappa$ can be adjusted so that a specific dimensionless Young's modulus $\tilde{E}$ is obtained that leads to small $|C_i|$. Changing this ratio is equivalent to moving along horizontal lines in the contour plots of Figure \ref{fig:errorConstants}. As an example and assuming a fixed grid spacing $\Delta x$, moving to the left in the plot corresponds to choosing smaller time steps and/or increased damping of the problem.
%Note that $\Delta t$ and $1/\kappa$ are numerically equivalent parameters and keeping the ratio $\Delta t/\kappa$ constant leads to the same relaxation times as is evident from Eqs. \eqref{eq:scaling_example} or \eqref{eq:scalingE}.

The visualization of the constants appearing in $r_x^{(4)}$ and $r_y^{(4)}$ shows that for each Poisson's ratio there exist parameter intervals of the dimensionless Young's modulus for which their values become small, thus we expect the leading-order error contribution to the numerical solution to be small as well. Further note that for some constants in Figure \ref{fig:errorConstants} an increase of its magnitude is observed for decreasing $\tilde{E}$. This implies that in general stronger damping or equivalently smaller time steps does not necessarily improve the numerical accuracy of the steady-state solution. This can be traced back to the fact that a different time step size changes the relaxation rates of the second-order moments, but leaves all other relaxation rates unaffected. In combination with the constant grid spacing this also violates the diffusive scaling on which the convergence property of the numerical scheme hinges. Similar observations in the context of fluid dynamics with the acoustic scaling and using the MRT collision operator have already been made in \citep{Dellar2003}.

A final simplification step assumes that all fourth-order derivatives of the solution components $u_x^{(0)}$ and $u_y^{(0)}$ are of similar magnitude. Under this assumption, a combined advantageous ratio of the discretization parameters $\Delta x$, $\Delta t$ and $\kappa$ can be estimated by minimizing the root sum squared of the constants $C_i$. If the influence of the leading-order error due to the body load is to be considered as well, the constants $D_i$ are also added. This in turn relies on the assumption that the spatial second-order derivatives of the forcing term are of similar magnitude as the fourth-order derivatives of the solution. The resulting error estimates $R_1$ and $R_2$ are given by
\begin{equation}
    \label{eq:rss_estimate}
    R_1 = \sqrt{\sum_{i=1}^5C_i^2} \qquad\qquad R_2 = \sqrt{\sum_{i=1}^5C_i^2+\sum_{i=1}^3D_i^2}.
\end{equation}

\begin{figure}[htbp]
    \centering
    \includegraphics{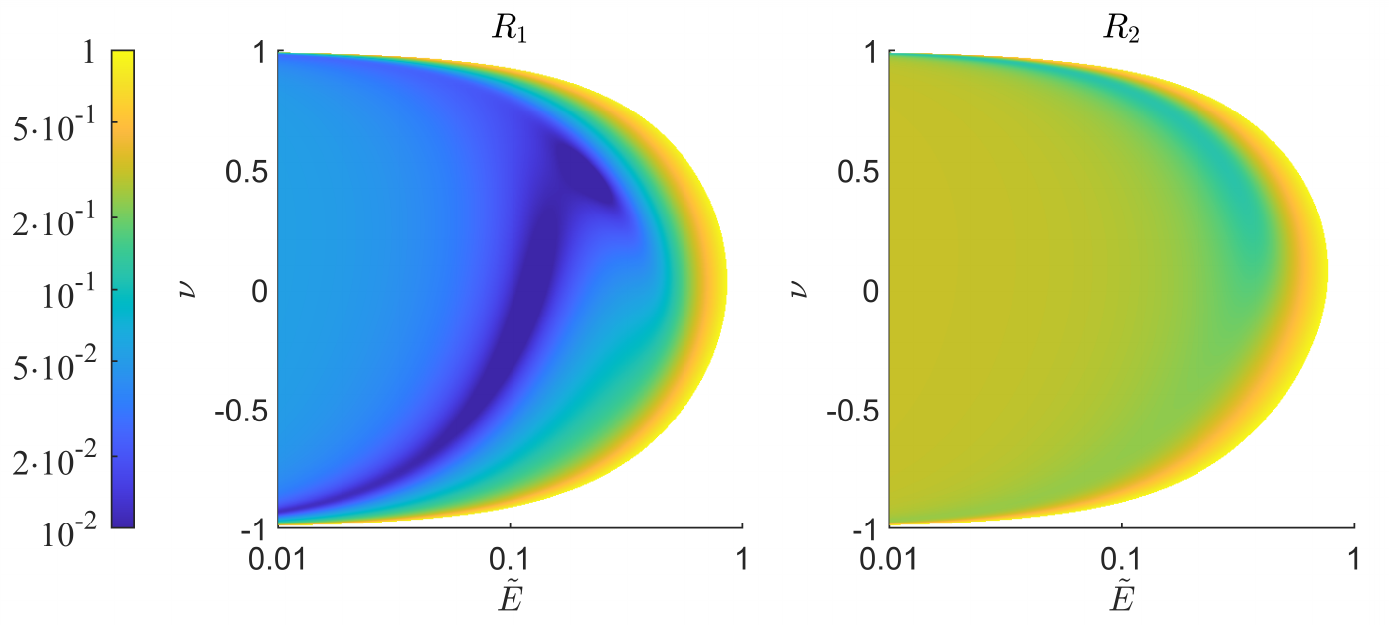}
    \caption[Leading-order error estimate]{Contour plot of the leading-order error estimates $R_1$ and $R_2$ as defined in \eqref{eq:rss_estimate}. Values above 1 cut off.}
    \label{fig:leadingError}
\end{figure}

Figure \ref{fig:leadingError} visualizes the influence of the dimensionless material parameters $\tilde{E}$ and $\nu$ on the qualitative estimate of the leading-order error. The direct comparison of the two contour plots showing $R_1$ and $R_2$ indicates that an error reduction of approximately one order of magnitude can be achieved by the partial removal of the leading-order error due to the body force in $r_x^{(4)}$ and $r_y^{(4)}$. Numerical experiments carried out in Section \ref{sec:num} will confirm this prediction.

Once again, the visualization of the combined constants shows that for each value of the Poisson's ratio there appears to be an advantageous interval of values for the dimensionless Young's modulus, which can be obtained by properly adjusting the relation between the discretization parameters in Eq. \eqref{eq:scalingE}.

For nearly incompressible material behavior -- corresponding to $\nu \rightarrow 1$ in 2D -- the optimal range of Young's moduli becomes relatively narrow and the smallest achievable value for $R_1$ increases. This in turn leads to larger expansion coefficients $u_x^{(2)}$ and $u_y^{(2)}$. Recalling the two possibilities of reducing the numerical error of the method, this observation shows that only a slight improvement following the strategy of the second option can be achieved. As a result, the numerical error of the method needs to be primarily reduced by the first option. Because this implies using a finer discretization, this can also be interpreted as decreasing numerical efficiency of the method for increasing Poisson's ratios. Numerical experiments in Section \ref{sec:num} will confirm this prediction as well.

\section{Linear stability analysis} \label{sec:stability}

As follows, we investigate the stability properties of the novel lattice Boltzmann scheme. The stability analysis follows hereby similar steps as in previous studies that all perform a linear von Neumann stability analysis  \citep{Sterling1996,Siebert2008,Hosseini2019}. In contrast to these works that study the stability properties of various lattice Boltzmann algorithms for fluid mechanics, no linearization about some homogeneous reference state is required because the present method is already linear. In brief, the stability analysis tests the attenuation or amplification properties for a given set of monochromatic planar waves by numerically computing the spectral radius of the linear operator that describes the lattice Boltzmann scheme \citep{Coreixas2020}. If the spectral radius becomes larger than 1, this implies that some modes of the distribution function grow for each iteration step and will eventually lead to an unbounded numerical solution, i.\,e. instability.

For the stability analysis, the populations are grouped in a vector-valued quantity $\boldsymbol{f}$, whose components $f_i,i=1,...,9$ correspond to the  populations $f_{kl}, k,l\in\{-1,0,1\}$ used so far. The mapping between the two-index notation used so far and the new single-index notation is described in Table \ref{tab:velocities} and follows the convention used e.\,g. in \citep{He1997} for the D2Q9 stencil. The same one-index notation is used in the following for the microscopic velocities $\boldsymbol{c_i}$. 
Note that, for the present scheme, the rest population $f_{00}$ can be removed as it has no influence on the behavior of the method. Accordingly, the stability analysis is performed on the resulting D2Q8 stencil.

Similarly, we introduce the vector $\boldsymbol{m}$ of the moments participating in the collision, whose components $m_i,i=1,...,8$ correspond to the moments $m_{\alpha}$ used so far as given in Table \ref{tab:moments}.

%\begin{table}[htbp]
%    \centering
%    \caption{Map from distribution index to streaming direction}
%    \begin{tabular}{r| *{8}{r}}
%        $i$ & 1 & 2 & 3 & 4 & 5 & 6 & 7 & 8 \\ \hline
%        $\boldsymbol{c}_{kl}$ & $\boldsymbol{c}_{10}$ & $\boldsymbol{c}_{01}$ & %$\boldsymbol{c}_{\bar{1}0}$ & $\boldsymbol{c}_{0\bar{1}}$ & $\boldsymbol{c}_{11}$ & %$\boldsymbol{c}_{\bar{1}1}$ & $\boldsymbol{c}_{\bar{1}\bar{1}}$ & %$\boldsymbol{c}_{1\bar{1}}$
%    \end{tabular}
%    \label{tab:velocities}
%\end{table}

\begin{table}[htbp]
    \centering
    \caption{Map between one-index and two-index notations for the populations}
    \begin{tabular}{r| *{8}{r}}
        $i$ & 1 & 2 & 3 & 4 & 5 & 6 & 7 & 8 \\ \hline
        $kl$ & $10$ & $01$ & $\bar{1}0$ & $0\bar{1}$ & $11$ & $\bar{1}1$ & $\bar{1}\bar{1}$ & $1\bar{1}$
    \end{tabular}
    \label{tab:velocities}
\end{table}

\begin{table}[htbp]
    \centering
    \caption{Map between one-index and two-index notations for the moments involved in the collision}
    \begin{tabular}{r|*{8}{r}}
        $i$ & 1 & 2 & 3 & 4 & 5 & 6 & 7 & 8 \\ \hline
        $\alpha$ & 10 & 01 & 11 & s & d & 12 & 21 & 22
    \end{tabular}
    \label{tab:moments}
\end{table}

The linear transformation between the vectors $\boldsymbol{f}$ and $\boldsymbol{m}$ is realized by the matrix $\boldsymbol{M}$ with components 
\begin{equation}
    \{M_{ij}\} = \left[\begin{array}{*{8}r}
    1 & 0 & -1 & 0 & \:\:1 & -1 & -1 & 1 \\
    0 & 1 & 0 & -1 & 1 & 1 & -1 & -1 \\
    0 & 0 & 0 & 0 & 1 & -1 & 1 & -1 \\
    1 & 1 & 1 & 1 & 2 & 2 & 2 & 2 \\
    1 & -1 & 1 & -1 & 0 & 0 & 0 & 0 \\
    0 & 0 & 0 & 0 & 1 & -1 & -1 & 1 \\
    0 & 0 & 0 & 0 & 1 & 1 & -1 & -1 \\
    0 & 0 & 0 & 0 & 1 & 1 & 1 & 1
    \end{array}\right],
\end{equation}
i.\,e. $m_i = \sum_{j=1}^8M_{ij}f_j$. $\boldsymbol{M}$ is easily constructed using the moment definition of Eq. \eqref{eq:moments} along with Eq. \eqref{eq:replace_raw}. 
In moment space, the collision involves the diagonal matrix $\boldsymbol{\Lambda} = \text{diag}(\omega_\alpha), \forall \alpha\in \{10,01,11,\text{s},\text{d},12,21,22\}$ with $\omega_{10}=\omega_{01}=0$. Naturally, the transformation back to the post-collision populations is carried out using the inverse of the transformation matrix, i.\,e. $\boldsymbol{M}^{-1}$. Lastly, the equilibrium moments are constructed from the populations by applying to the population vector the matrix $\boldsymbol{M}^{eq}$ with components
\begin{equation}
    \{M_{ij}^{eq}\} = \left[\begin{array}{*{8}r}
    1 & 0 & -1 & 0 & \:\:1 & -1 & -1 & 1 \\
    0 & 1 & 0 & -1 & 1 & 1 & -1 & -1 \\
    0 & 0 & 0 & 0 & 0 & 0 & 0 & 0 \\
    0 & 0 & 0 & 0 & 0 & 0 & 0 & 0 \\
    0 & 0 & 0 & 0 & 0 & 0 & 0 & 0 \\
    0 & 0 & 0 & 0 & \theta & -\theta & -\theta & \theta \\
    0 & 0 & 0 & 0 & \theta & \theta & -\theta & -\theta \\
    0 & 0 & 0 & 0 & 0 & 0 & 0 & 0
    \end{array}\right].
\end{equation}

Using these definitions, one iteration step of the lattice Boltzmann method can be described through a linear operator as follows
\begin{equation}
    f_i(\tilde{\boldsymbol{x}}+\tilde{\boldsymbol{c}}_i,\tilde{t}+1) = \sum_{j,k,l=1}^8\underbrace{\left[M^{-1}_{ij}\Lambda_{jk}M_{kl}^{eq} + M^{-1}_{ij}(\delta_{jk} - \Lambda_{jk})M_{kl}\right]}_{\textstyle\mathcal{A}_{il}} f_l (\tilde{\boldsymbol{x}},\tilde{t}),
    \label{eq:lbmvector}
\end{equation}
where the linear collision operator $\mathcal{A}_{ij}$ has been introduced. The properties of the corresponding matrix can be investigated using linear algebra, giving immediate information on the stability of the method. As in \citep{Hosseini2019,Coreixas2020}, the amplification or attenuation properties of the lattice Boltzmann scheme are investigated by assuming that each population in $\boldsymbol{f}$ has the form of a monochromatic planar wave with wave vector $\boldsymbol{k} = [k_x\;k_y]^T$, complex frequency $\varpi$ and amplitude $a_i$ 
\begin{equation}
    f_i(\tilde{\boldsymbol{x}},\tilde{t}) = a_i \exp{\left[\imath(\boldsymbol{k}\cdot\tilde{\boldsymbol{x}}-\varpi \tilde{t})\right]} \quad \Leftrightarrow \quad  \boldsymbol{f}(\tilde{\boldsymbol{x}},\tilde{t}) = \boldsymbol{a} \exp{\left[\imath(\boldsymbol{k}\cdot\tilde{\boldsymbol{x}}-\varpi \tilde{t})\right]},
\end{equation}
where $\imath$ is the imaginary unit. Note that for the stability analysis $c=\Delta x=\Delta t = 1$ is assumed without loss of generality, because this can be interpreted as a re-scaling of the wave vector $\boldsymbol{k}$ and the complex frequency $\varpi$. Introducing this ansatz into Eq. \eqref{eq:lbmvector} leads to the following eigenvalue problem
\begin{equation}
    \boldsymbol{L} \boldsymbol{a} = \exp{(-\imath \varpi)} \boldsymbol{a},
    \label{eq:eigen}
\end{equation}
with the matrix $\boldsymbol{L}$ defined as
\begin{equation}
    L_{ij} = \exp{(-\imath \boldsymbol{k} \cdot \tilde{\boldsymbol{c}}_i)}\mathcal{A}_{ij}.
\end{equation}

Instability is encountered if some mode of $\boldsymbol{f}$ is amplified during each iteration. Accordingly, the following condition needs to hold for a stable method:
\begin{equation}
    \sup_{\boldsymbol{k}}\rho(\boldsymbol{L}) \leq 1,
\end{equation}
where $\rho(\boldsymbol{L})$ denotes the spectral radius of $\boldsymbol{L}$.

The condition above is approximated by sampling through the set of admissible wave vectors and numerically computing the largest eigenvalue of $\boldsymbol{L}$. For this purpose, the following parameterization is adopted for the wave vectors in 2D
\begin{equation}
    \boldsymbol{k} = \hat{k}\left[\begin{matrix}
    \cos{\varphi} \\ \sin{\varphi}
    \end{matrix}\right] \qquad \text{with } \hat{k} \in [0,\pi] \quad \text{and } \varphi \in [0,\pi/4],
\end{equation}
where $\hat{k}$ is the magnitude of the wave vector and $\varphi$ the angle between the wave vector and the x-coordinate direction. The ranges of the two parameters follow from symmetry considerations \citep{Siebert2008,Hosseini2019}. In the long wave length limit, i.\,e. $\boldsymbol{k}\rightarrow\boldsymbol{0}$, the classical result for the relaxation rates is retrieved in order to guarantee stability, i.\,e.
\begin{equation}
    0 \leq \omega_\alpha \leq 2 \qquad \forall \alpha \in \mathcal{I}.
    \label{eq:stability}
\end{equation}

For the case of arbitrary wave vectors, it has been observed that, for the present method, instability is always first encountered for axis-aligned wave vectors, i.\,e. $\varphi=0$. For this reason, only 5 points are sampled from the domain of definition of $\varphi$. In order to sufficiently resolve the space of possible wave vectors the magnitude $\hat{k}$ is evaluated at 50 points in $[0,\pi]$. This procedure is repeated for a grid of values covering the parameter space spanned by the dimensionless Young's modulus and the Poisson's ratio. For all other parameters of the collision operator, the standard values as in Table \ref{tab:scheme} are set.

\begin{figure}[htbp]
    \centering
    \begin{minipage}[t]{0.5\textwidth}
    \begin{flushright}
        \vspace{0pt}
        \includegraphics{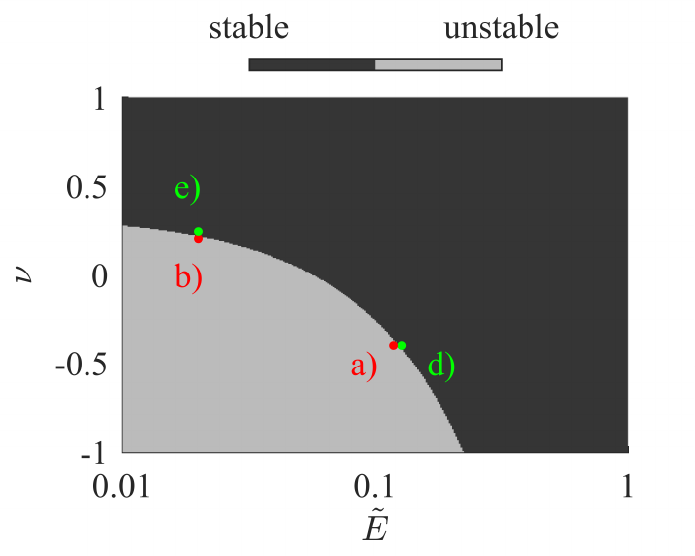}
    \end{flushright}
    \end{minipage}
    \begin{minipage}[t]{0.45\textwidth}
        \begin{flushleft}
        \vspace{21.7pt}
        \includegraphics{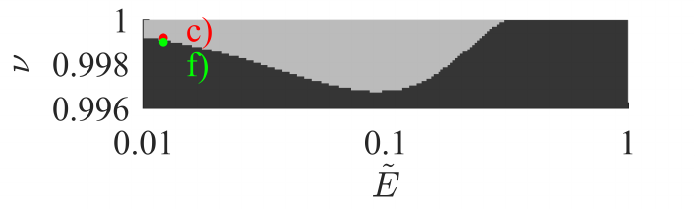}
    \end{flushleft}
    \end{minipage}
    \caption[Stable and unstable regions in material parameter space]{Stable and unstable regions in material parameter space. Left: Complete range of physically admissible values of the Poisson's ratio. Right: Inset for nearly incompressible material behavior. The dots indicate material parameters used for numerical verification in Section \ref{sec:num}.}
    \label{fig:stability}
\end{figure}

The left side of Figure \ref{fig:stability} shows the stability region for all admissible Poisson's ratios and the range of values for the dimensionless Young's modulus that is relevant for the simulations. The analysis indicates that there exists a region with unstable behavior for small values of the dimensionless Young's modulus and negative or small positive Poisson's ratio. Recalling the leading-order error analysis of the previous section it becomes apparent that the onset of instability poses a practical limit only for sufficiently small Poisson's ratios with approximately $\nu < -0.3$. For larger Poisson's ratios the predicted beneficial combination of discretization parameters as demonstrated in Figure \ref{fig:leadingError} lies in the stable region. However, for cases with $\nu < -0.3$ the onset of instability inhibits the use of a beneficial discretization to reduce the leading-order error influence. As a result, problems with $\nu < -0.3$ can only be solved with reduced accuracy or require significantly finer meshes to achieve the same accuracy. However, negative values of $\nu$, while thermodynamically possible, are irrelevant for most practical purposes.

The right-hand side of Figure \ref{fig:stability} shows another region of instability encountered for nearly incompressible material behavior. However, the instability occurs only approximately for $\nu > 0.996$ for which the leading-order error analysis predicts very large errors. Therefore, this regions poses no practical limitations on the range of Poisson's ratios the method can handle. Altogether, the combined results of the stability analysis and the leading-order error investigation predict that the novel method is accurate and stable for approximately $\nu \in [-0.3,0.95]$.

\section{Periodic boundary and initial conditions} \label{sec:icbc}

The application of physically consistent, accurate and stable boundary conditions for the lattice Boltzmann method is a challenging task, because the physical conditions need to be translated into expressions for the populations. For this reason, the influence of boundary conditions is not considered in this initial contribution and only problems with periodic boundary conditions are solved. To enforce periodicity, a simple formulation from the literature can be directly applied as will be outlined in the following.

Because the problem is solved by advancing in pseudo-time, an initial condition for the populations needs to be specified as well. For this purpose, the popular initialization at local equilibrium is investigated using the asymptotic expansion of Section \ref{sec:asympt}.

\subsection{Periodic boundary}

For the periodic problem, axis-aligned rectangular domains of size $L_x \times L_y$ are considered. The respective boundary conditions are realized by copying the outgoing populations to the opposite end of the domain \citep{Kruger2017}. In order to avoid duplicating the nodes on the boundaries that are connected by periodicity, the node lattice is offset by half a grid spacing from the physical domain boundary.

%(see Figure \ref{fig:periodicity} for an example involving the left and right boundaries). For an efficient implementation it is beneficial to keep the streaming step unchanged throughout the whole computational domain. In order to achieve this, a single layer of ghost nodes is appended to all sides of the domain into which the outgoing distributions are streamed. A subsequent step performed only in the ghost layers executes the transport of the outgoing distributions to the opposite end of the domain \citep{Kruger2017}.
%\begin{figure}[htbp]
%    \centering
%    \includegraphics{figures/periodicity.pdf}
%    \caption{Periodic boundary condition shown for the left-right side pair}
%    \label{fig:periodicity}
%\end{figure}

\subsection{Initial condition}

A popular method to enforce the initial condition $\tilde{\boldsymbol{u}}=\tilde{\boldsymbol{u}}_0$ (Eq. \eqref{eq:initnondim}) is to initialize all populations $f_{ij}$ at equilibrium. This is equivalent to prescribing the local equilibrium value for all moments (see Table \ref{tab:scheme}) followed by a back-transformation into populations as done during the collision stage. The respective initialization of all moments is given by \eqref{eq:init_cond}.
\begin{equation}
    \label{eq:init_cond}
    m_\alpha(\tilde{\boldsymbol{x}},0) = m^{eq}_\alpha(\tilde{\boldsymbol{u}}_0(\tilde{\boldsymbol{x}})) \qquad \forall\alpha\in \{10,01,11,\text{s},\text{d},12,21,22\}
\end{equation}

In order to investigate the consistency order of Eq. \eqref{eq:init_cond}, the asymptotic expansion of all moments is computed up to first order using Eq. \eqref{eq:expansion} and the results in Section \ref{sec:asympt}.
\begin{alignat}{6}
\label{eq:pre_moments_first}
m_{10} &= u_x^{(0)} &&+ \varepsilon u_x^{(1)} &&+ \mathcal{O}(\varepsilon^2) \\
m_{01} &= u_y^{(0)} &&+ \varepsilon u_y^{(1)} &&+ \mathcal{O}(\varepsilon^2) \\
m_{11} &=  &&- \varepsilon \omega_{11}^{-1}\theta\left(\partial_{\tilde{x}}u_y^{(0)}+\partial_{\tilde{y}}u_x^{(0)}\right) &&+ \mathcal{O}(\varepsilon^2) \\
m_\text{s} &=  &&- \varepsilon \omega_\text{s}^{-1}(1+\theta)\left(\partial_{\tilde{x}}u_x^{(0)}+\partial_{\tilde{y}}u_y^{(0)}\right) &&+ \mathcal{O}(\varepsilon^2) \\
m_\text{d} &=  &&- \varepsilon \omega_\text{d}^{-1}(1-\theta)\left(\partial_{\tilde{x}}u_x^{(0)}-\partial_{\tilde{y}}u_y^{(0)}\right) &&+ \mathcal{O}(\varepsilon^2) \\
m_{12} &= \theta u_x^{(0)} &&+ \varepsilon \theta u_x^{(1)} &&+ \mathcal{O}(\varepsilon^2) \\
m_{21} &= \theta u_y^{(0)} &&+ \varepsilon \theta u_y^{(1)} &&+ \mathcal{O}(\varepsilon^2) \\
m_{22} &=  &&- \varepsilon \omega_{22}^{-1}\theta\left(\partial_{\tilde{x}}u_x^{(0)}+\partial_{\tilde{y}}u_y^{(0)}\right) &&+ \mathcal{O}(\varepsilon^2)
\label{eq:pre_moments_last}
\end{alignat}
Introducing these results into Eq. \eqref{eq:init_cond} reveals that second-order consistency with $\boldsymbol{u}^{(1)}=\boldsymbol{0}$ can only be achieved with this initialization if $\nabla \boldsymbol{\tilde{u}}_0 = \boldsymbol{0}$. Because all numerical examples in Section \ref{sec:num} start from $\tilde{\boldsymbol{u}}_0=\boldsymbol{0}$, this simplified initial condition is sufficient.

\section{Numerical verification} \label{sec:num}

This section presents a few numerical examples. These serve the purpose of verifying the analytical derivations of the previous sections and thus of assessing the performance of the new scheme. To measure the error in the displacement and the stress solution, some grid norm definitions are introduced in the first subsection. In order to allow for as general conclusions as possible, the accuracy of the method is investigated using manufactured solutions. This concept is briefly introduced in the following subsection. Using the method of manufactured solutions, a series of convergence studies are carried out with the aim of verifying the results of the asymptotic expansion. By considering only periodic problems, the findings of the leading-order error and stability analysis of Section \ref{sec:conv} can be numerically explored in great detail.

\subsection{Grid norm definitions}

Before defining the grid norms, the numerical (discretization) error in the displacement and the stress solution is introduced. As outlined in Section \ref{sec:asympt}, the approximate displacement solution is obtained from the first-order moments, which are defined on the set of all lattice node positions $\{\boldsymbol{x}_i\,|\,\boldsymbol{x}_i\in\Omega\}_{i=1}^N$. $N$ is the total number of grid nodes contained in the problem domain $\Omega$. Keeping in mind that the present method advances in pseudo-time steps $t_j \in [0,t_f]$ until steady state is achieved within a specified tolerance, the numerical solution also depends on time. Provided that the exact solution is available at each lattice node $\hat{\boldsymbol{u}}_i = \hat{\boldsymbol{u}}(\boldsymbol{x}_i)$, the error in the displacement field solution is defined as
\begin{equation}
    \boldsymbol{e}^{\boldsymbol{u}}_i(t_j) = \boldsymbol{u}^{num}_i(t_j) - \hat{\boldsymbol{u}}_i
    \qquad \text{with } \boldsymbol{u}^{num}_i = U\left[ \begin{matrix}
    \bar{m}_{10} \\ \bar{m}_{01}
    \end{matrix} \right](\boldsymbol{x}_i,t_j).
\end{equation}
Note that the numerical approximation of the solution, obtained in dimensionless form, is transformed back into physical units for the error evaluation.

In order to simplify the grid norm computation of the error in the stress solution, the independent stress components are regrouped into a vector. As derived in section \ref{sec:stress}, the components of the approximate Cauchy stress are identified with the bared second-order  moments, thus 
\begin{equation}
    \boldsymbol{e}^{\boldsymbol{\sigma}}_i(t_j) = \boldsymbol{\sigma}^{num}_i(t_j) - \hat{\boldsymbol{\sigma}}_i
    \qquad \text{with } \boldsymbol{\sigma}^{num}_i(t_j) = -LT^{-1}U\kappa\left[ \begin{matrix}
    (\bar{m}_\text{s} + \bar{m}_\text{d})/2 \\ (\bar{m}_\text{s} - \bar{m}_\text{d})/2 \\ \bar{m}_{11}
    \end{matrix} \right](\boldsymbol{x}_i,t_j).
\end{equation}
Once again, the numerical result is transformed into physical units (see Section \ref{sec:nondim} for the scaling factors). Here $\hat{\boldsymbol{\sigma}}_i$ denotes the exact stress solution at lattice node $i$.
As a next step, two grid norms are introduced to obtain global measures of the previously defined error functions. The norm definitions derive from the continuous $L_2(\Omega)$ and $L_\infty(\Omega)$ norms by applying numerical quadrature at the grid nodes. Because the norm computation involves a summation only over all nodes in space, the error measure remains a function of time that indicates the convergence of the method towards the static equilibrium solution. Using the example of the displacement error, the two norm definitions are provided below
\begin{alignat}{5}
    &\text{L2}\qquad &&\lVert \boldsymbol{e}^{\boldsymbol{u}} \rVert_2(t_j) &&= \left( (\Delta x)^2 \sum_{i=1}^N |\boldsymbol{e}^{\boldsymbol{u}}_i|^2_2(t_j) \right)^{\frac{1}{2}} \\
    &\text{Linf}\qquad &&\lVert \boldsymbol{e}^{\boldsymbol{u}} \rVert_\infty(t_j) &&= \max_{i=1,\dots,N} |\boldsymbol{e}^{\boldsymbol{u}}_i|_\infty(t_j),
\end{alignat}
where $|\cdot|_2$ and $|\cdot|_\infty$ denote the Euclidean norm and the maximum norm of a vector. The following observations apply:
\begin{itemize}
\item $\Delta x$ is the uniform grid spacing of the lattice and scales proportionally with the smallness parameter of the asymptotic expansion in Section \ref{sec:asympt}, i.\,e. $\Delta x \sim \varepsilon$.
\item Because the scheme is designed to determine the solution at static equilibrium (Eq. \eqref{eq:linelas}), the grid norm of the error is considered at final time $t_f$ unless stated otherwise. %i.\,e. $\lVert \boldsymbol{\varphi} \rVert_{p} = \lVert \boldsymbol{\varphi} \rVert_{p}(t_f),p\in\{2,\infty\}$.
%\item[Remark 4] The Linf norm is a very sensitive measure of the error, because it is governed by the largest component of the vector-valued error at the grid node with the maximum deviation. Thus, establishing the predicted convergence order in this norm can be considered a strong confirmation of the theory. The L2 norm on the other hand provides a measure of the error in an average-sense, which helps to estimate the expected error of the method for similar problems.
\item In order to provide easily interpretable results, in the following we report the relative error, i.\,e. the absolute error divided by the L2 grid norm of the exact solution, e.\,g.: $\lVert \boldsymbol{e}^{\boldsymbol{u}} \rVert_\infty / \lVert \hat{\boldsymbol{u}} \rVert_2$.
\end{itemize}

\subsection{Method of manufactured solutions}

As follows, we briefly outline the method of manufactured solutions \citep{Roache2002} that can generate almost arbitrary exact solutions to be used for the numerical verification. As a starting point for the method, the desired exact solution $\hat{\boldsymbol{u}}$ is freely chosen, but needs to be sufficiently differentiable so that it can be a solution to the problem. The next step involves computing the matching source term $\boldsymbol{b}$ so that $\hat{\boldsymbol{u}}$ actually solves the governing equation in conjunction with the source term. In the present context (see \eqref{eq:linelasdamp}) this implies computing a body load as shown below:
\begin{equation}
\label{eq:manu_eq}
    \boldsymbol{b} = \kappa \partial_t \hat{\boldsymbol{u}} - \mu \nabla^2 \hat{\boldsymbol{u}} - K \nabla (\nabla \cdot \hat{\boldsymbol{u}}).
\end{equation}
For the case of quasi-static linear elasticity, time-independent solutions are assumed so that the first term on the right-hand side of Eq. \eqref{eq:manu_eq} vanishes.

Within the context of periodic problems, the manufactured solution ansatz needs to satisfy two additional constraints. The first condition is given by Eq. \eqref{eq:periodicity}. The second requirement comes from the fact that the solution of the periodic problem conserves the mean value of the initialization. Because all examples in this section are initialized with $\boldsymbol{u}_0 = \boldsymbol{0}$, the manufactured solution needs to have zero mean.

With the body load and initial condition defined, the numerical solution is computed and can be compared against the exact analytical solution using the grid norm definitions introduced before.

\subsection{Numerical examples}

In this section, a selection of numerical examples is presented that solve the periodic problem of Eqs. \eqref{eq:linnondim}--\eqref{eq:periodicitynondim} for a given body load generated through the method of manufactured solutions. The first example demonstrates the second-order convergence property of the method in both the L2 and Linf norm for the displacement and the stress solution. The second example showcases the improvement of the accuracy of the method if the leading-order error contribution due to the body load is compensated for. This is followed by an example that is designed in such a way that fourth-order consistency is numerically observed if the analytical derivations in Section \ref{sec:conv} are correct. Furthermore, the effectiveness of the error estimate to identify a beneficial ratio of grid spacing and time step size is investigated by numerical parameter studies that probe the actual error of the method for various combinations of discretization parameters. Subsequently, the analytically predicted location of the stability boundary in material parameter space is verified for a few testing points and finally the results obtained with the novel lattice Boltzmann scheme are compared with linear finite element results.

\subsubsection{Convergence study of the standard scheme}

\begin{figure}[htbp]
    \centering
    \includegraphics{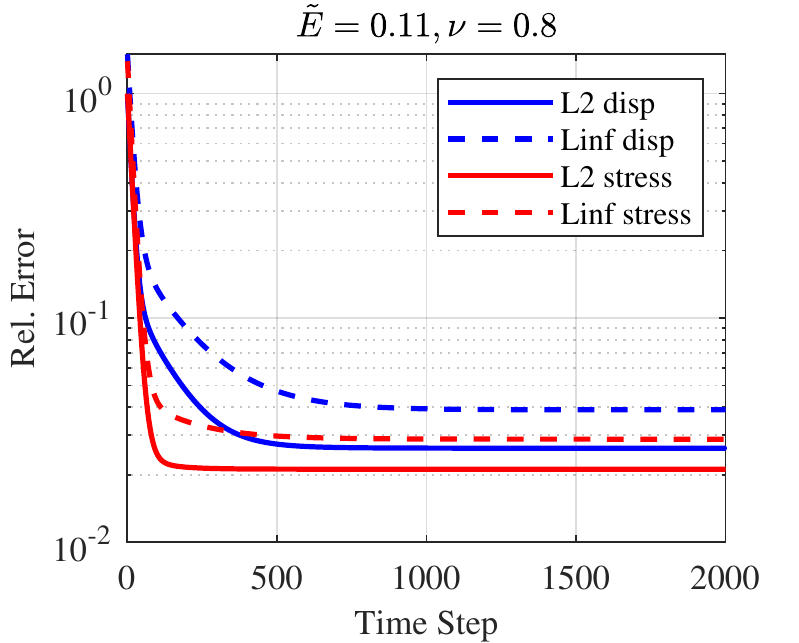}
    \caption{Temporal evolution of the error}
    \label{fig:evolution_periodic}
\end{figure}

As a first example, the periodic problem is solved with the standard numerical scheme, i.\,e. the scheme that does not compensate for the leading-order error due to the body load (see Section \ref{sec:conv}). In this case the material parameters are chosen as $\nu=0.8$ and $\tilde{E}=0.11$ to obtain a small leading-order error contribution based on the investigation of Section \ref{sec:conv}. The manufactured solution on a $[0,L]\times[0,L]$ domain used for the following convergence study is
\begin{alignat}{6}
    \label{eq:x_manu_standard}
    \hat{u}_x/U &= 9\cdot 10^{-4} &&\cos{(2\pi \tilde{x})} &&\sin{(2\pi \tilde{y})} \\
    \label{eq:y_manu_standard}
    \hat{u}_y/U &= 7\cdot 10^{-4} &&\sin{(2\pi \tilde{x})} &&\cos{(2\pi \tilde{y})},
\end{alignat}
with $\tilde{x}=x/L$ and $\tilde{y}=y/L$ as introduced in Section \ref{sec:nondim}. The final time $t_f$ is set sufficiently large so that steady state is reached. An example of the convergence of the numerical solution towards steady state over the pseudo-time steps is provided by Figure \ref{fig:evolution_periodic}. It shows the evolution of the relative error of the displacement and the stress solution in the L2 and the Linf norm. Note that this simulation is performed on the coarsest mesh of the following convergence study with $\Delta x = 0.05L$, for which relative errors of approximately $3 - 4\%$ in the Linf norm are obtained for the given discretization.

Steady state is reached after roughly $10^3$ time steps and the stress field converges faster than the displacement solution. Despite the large number of time steps, the total computing time of this example is still very low: ca. $1.6$\,s on a low-end notebook with a non-optimized implementation. 

\begin{figure}[htbp]
    \centering
    \begin{minipage}{0.5\textwidth}
    \centering
    \includegraphics{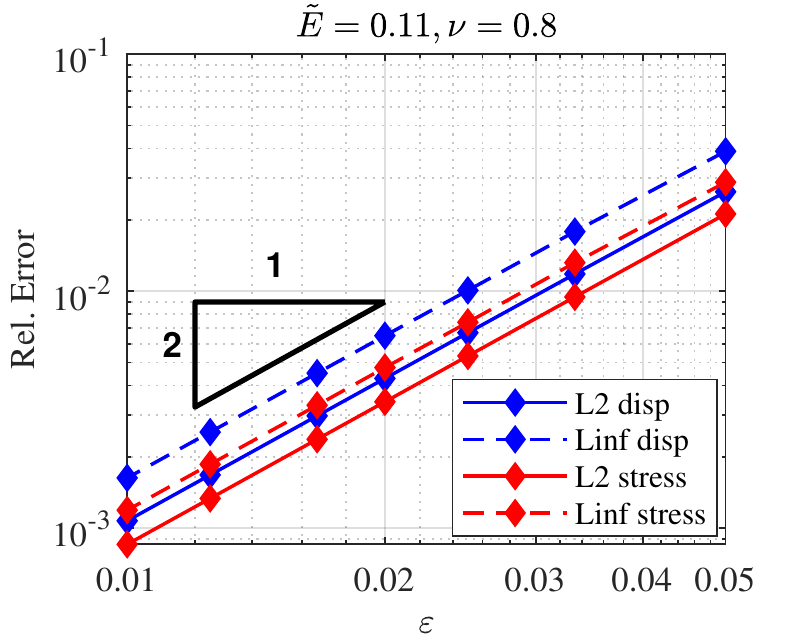}
    \end{minipage}
    \begin{minipage}{0.3\textwidth}
   % \begin{tabular}{l  l}
%    Norm & Conv. Order \\ \hline
%        L2 disp & 2.01 \\
%        Linf disp & 2.00 \\
%        L2 stress & 2.00 \\
%        Linf stress & 1.98
%    \end{tabular}
    \end{minipage}
    \caption{Convergence study for the implementation without compensation of the leading-order error due to the body load}
    \label{fig:periodic_standard}
\end{figure}

Using the already described problem definition, a grid convergence study is carried out. As predicted by the asymptotic analysis in Sections \ref{sec:asympt}, \ref{sec:stress} and \ref{sec:conv}, the results in Figure \ref{fig:periodic_standard} indicate second-order convergence for both the displacement and the stress solution for decreasing smallness parameter $\varepsilon=\Delta x/L$. For these results the total number of grid nodes is varied between $20^2$ and $100^2$ nodes. Notably, a second-order convergence rate is also achieved in the Linf norm for both solution quantities, which confirms the analytical derivations with high accuracy.

\subsubsection{Convergence study of the improved scheme} \label{sec:improved_scheme}

\begin{figure}[htbp]
    \centering
    \begin{minipage}{0.5\textwidth}
    \centering
    \includegraphics{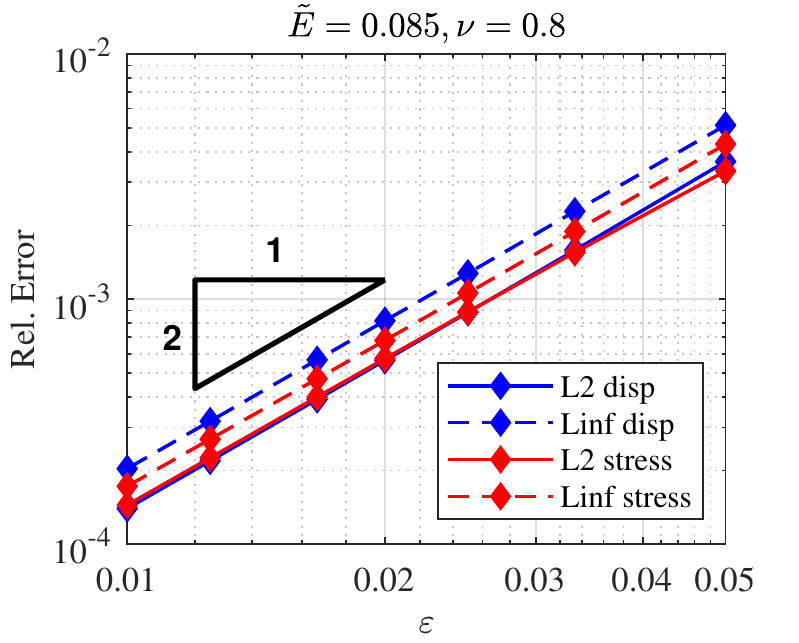}
    \end{minipage}
    \begin{minipage}{0.3\textwidth}
%    \begin{tabular}{l  l}
%    Norm & Conv. Order \\ \hline
%        L2 disp & 2.04 \\
%        Linf disp & 2.03 \\
%        L2 stress & 1.95 \\
%        Linf stress & 2.01
%    \end{tabular}
    \end{minipage}
    \caption{Convergence study for the implementation with compensation of leading-order error due to the body load}
    \label{fig:periodic_improved}
\end{figure}

As suggested in Section \ref{sec:conv}, the portion of the leading-order error due to the body load can be compensated for by a higher-order correction of the forcing term. With this modification, the optimal value of the dimensionless Young's modulus is predicted to be slightly different based on the method in Section \ref{sec:conv}. Therefore, the simulation is performed for $\nu=0.8$ and $\tilde{E}=0.085$. Besides this adjustment, the same problem as for the previous case is solved. Figure \ref{fig:periodic_improved} shows the results for the scheme with compensation of the leading-order error due to the body load. Once again, approximately second-order convergence is observed from the numerical experiments, because the leading-order error is only partially removed. However, comparing the results of the convergence studies for the two variants (Figures \ref{fig:periodic_standard} and \ref{fig:periodic_improved}) shows a reduction of the errors by roughly one order of magnitude. This improvement was already concluded from the derivations in Section \ref{sec:conv} and specifically by Figure \ref{fig:leadingError}. In summary, a significant improvement of the numerical accuracy can be expected for this partial compensation of the leading-order error that requires only minimal extra computational effort.

\subsubsection{A special case with fourth-order convergence}

The next convergence study on periodic domains aims at constructing a special case in which, according to the asymptotic expansion, fourth-order convergence is expected, and at verifying this expectation numerically. To this end, an analytical solution with vanishing mixed spatial derivatives is chosen as follows
\begin{alignat}{4}
    \hat{u}_x/U &= 9 \cdot 10^{-4} \left(\cos{(2\pi \tilde{x})}\right. &&+ \left.\sin{(2\pi \tilde{y})}\right) \\
    \hat{u}_y/U &= 7 \cdot 10^{-4} \left(\sin{(2\pi \tilde{x})}\right. &&+ \left.\cos{(2\pi \tilde{y})}\right).
\end{alignat}
As a result, only the two terms with the constants $C_1$ and $C_5$ influence the leading-order error of both $u_x^{(2)}$ and $u_y^{(2)}$ as can be seen from Eq. \eqref{eq:ux_order2} for the x-component. Therefore, only two requirements ($C_1=0$ and $C_5=0$) need to be satisfied to remove the $r_x^{(4)}$ and $r_y^{(4)}$ terms, which is less than the number of free parameters in the MRT collision operator (see Section \ref{sec:conv}) so that a complete removal of $r_x^{(4)}$ and $r_y^{(4)}$ is theoretically possible. For this example the relaxation times $\tau_{12}=\tau_{21}$ and $\tau_{22}$ have been appropriately chosen to achieve this and are reported in \ref{sec:rates_fourth_order}. Furthermore, the asymptotic analysis also shows that the next higher-order error contribution vanishes if $u_x^{(1)}=u_y^{(1)}=0$ (not demonstrated here). This is fulfilled by the second-order consistency of the method so that a fourth-order accurate method can be expected.

\begin{figure}[htbp]
    \centering
    \begin{minipage}{0.5\textwidth}
    \centering
    \includegraphics{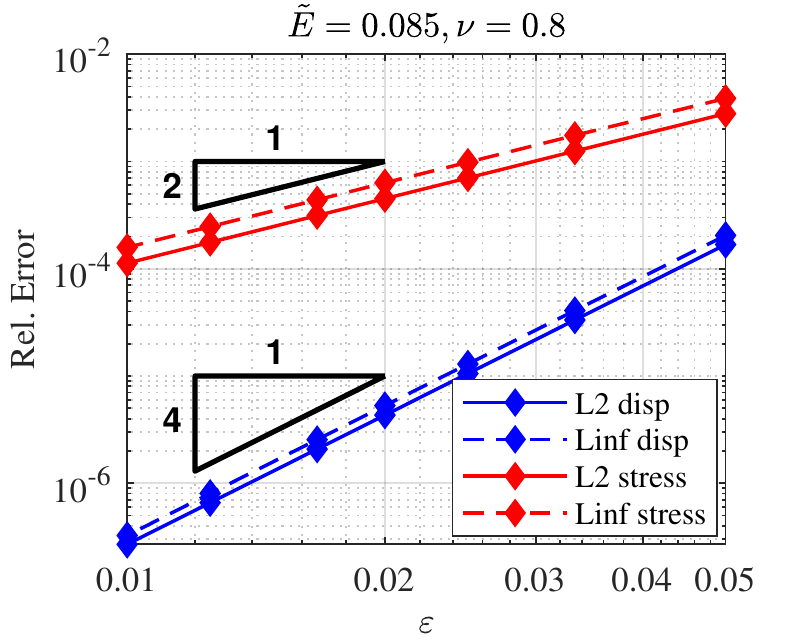}
    \end{minipage}
    \begin{minipage}{0.3\textwidth}
%    \begin{tabular}{l  l}
%    Norm & Conv. Order \\ \hline
%        L2 disp & 4.00 \\
%        Linf disp & 4.00 \\
%        L2 stress & 1.99 \\
%        Linf stress & 1.99
%    \end{tabular}
    \end{minipage}
    \caption{Convergence study for the special case with fourth-order consistency}
    \label{fig:periodic_fourth}
\end{figure}

Figure \ref{fig:periodic_fourth} shows the numerically obtained fourth-order convergence in the displacement field solution in both norms, which is a strong confirmation of the analytical derivations in Section \ref{sec:conv}. It is also observed that the leading-order error correction for the displacement solution does not significantly affect the accuracy of the numerical stress solution.

Unfortunately, this fourth-order consistent scheme requires a special structure of the exact solution, specifically that all mixed fourth-order spatial derivatives vanish. Furthermore, the higher-order relaxation rates required to let $C_1$ and $C_5$ vanish can be shown to violate the stability condition of Eq. \eqref{eq:stability} for a large region in the space of the admissible combinations of $\tilde{E}$ and $\nu$. Therefore, this higher-order accurate method is only applicable to a restricted class of problems.

\subsubsection{Leading-order error behavior for different combinations of discretization parameters}

During the convergence studies with the second-order consistent schemes, the magnitude of the leading-order error contribution was already reduced by appropriately choosing $\tilde{E}$ (which, for given $E$ and lattice spacing, amounts to choosing $\kappa$ or $\Delta t$) . The purpose of the following study is to investigate the accuracy and reliability of the leading-order error estimate $R_1$ that has been introduced in in Eq. \eqref{eq:rss_estimate} in Section \ref{sec:stability}. To this end, parameter studies are carried out for several values of the Poisson's ratio, by varying the time scaling in order to realize different values for the dimensionless Young's modulus. This is equivalent to moving along horizontal lines in the contour plot of Figure \ref{fig:leadingError}. Further note that all simulations are performed using the scheme with compensation of the leading-order error due to the body load. Using the manufactured solution of Eqs. \eqref{eq:x_manu_standard} and \eqref{eq:y_manu_standard}, the numerical error of the displacement solution is evaluated in the L2 norm for a range of different combinations of discretization parameters and compared with the error estimate. The value of the error estimate $R_1$ (see Eq. \eqref{eq:rss_estimate}) and the actual numerical error does not match in general, because the latter is (to leading-order) the solution to Eq. \eqref{eq:ux_order2}, whereas the former is computed from the coefficients appearing in the body load term of the same governing equation. However, it is expected that the qualitative shape of the error graphs agrees with each other. Specifically, the discretization for which the estimate predicts the smallest leading-order error contribution should lie in the vicinity of the value for which the actual numerical error has its minimum.

\begin{figure}[htbp]
    \centering
    \includegraphics{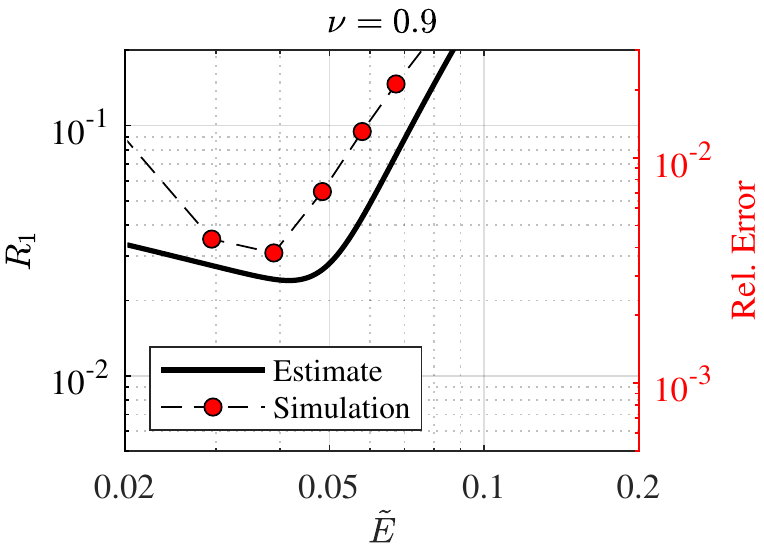}\;
    \includegraphics{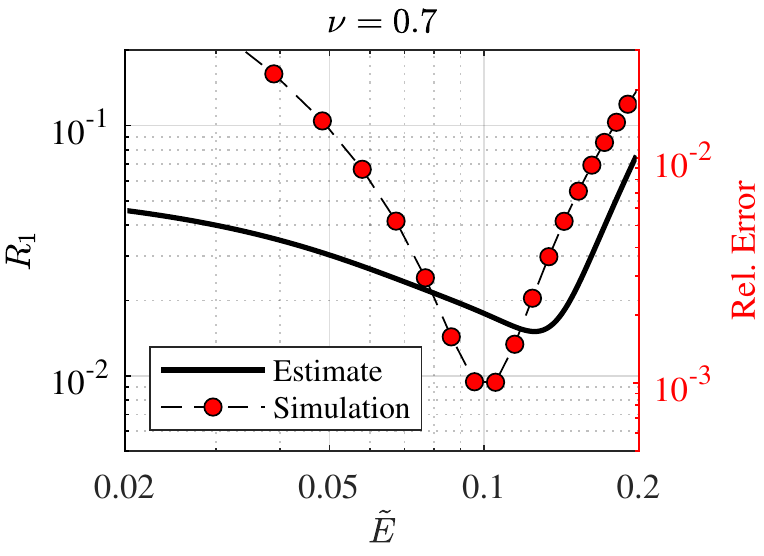}
    \includegraphics{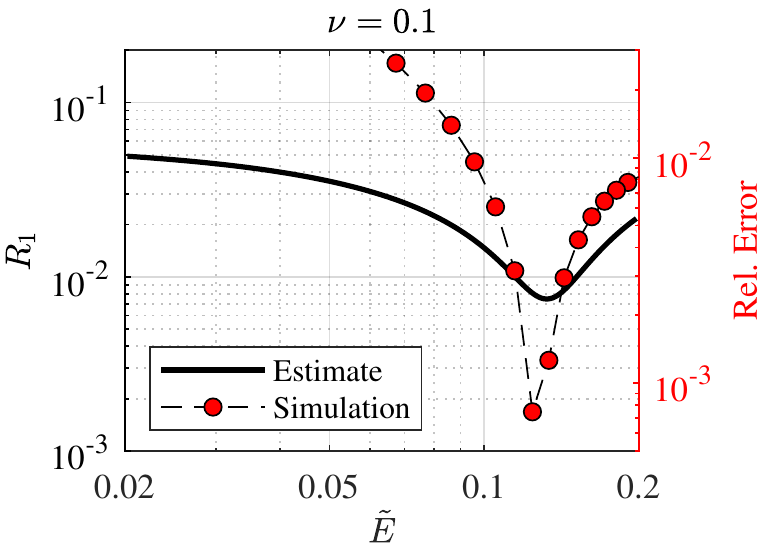}\;
    \includegraphics{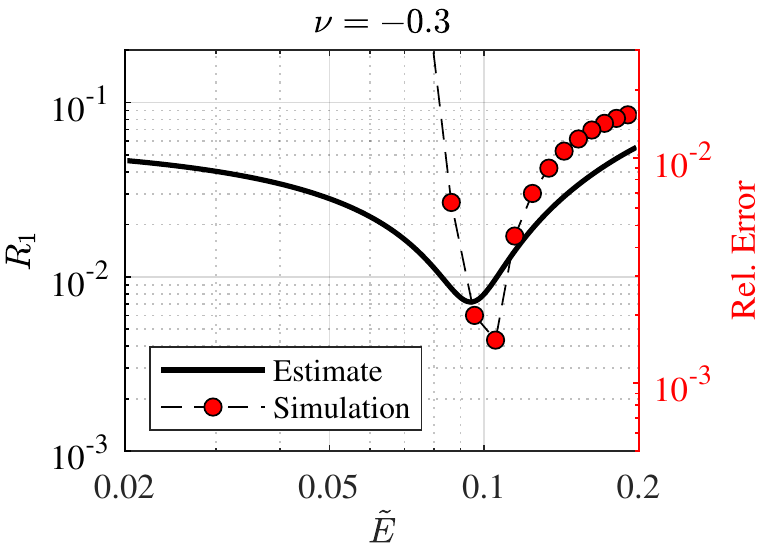}
    \caption{Qualitative comparison of the L2 error in the displacement solution (right y-axis) with leading-order error estimate (left y-axis)}
    \label{fig:num_leadingError}
\end{figure}

This comparison is presented in Figure \ref{fig:num_leadingError} at four different values of the Poisson's ratio. The examples demonstrate a satisfactory qualitative agreement between the error estimate $R_1$ (see Eq. \eqref{eq:rss_estimate}) and the actual error. Comparing the results for the different Poisson's ratios shows that the estimate predicts slightly too large values for the dimensionless Young's modulus in the case of large Poisson's ratios and vice versa. When comparing the smallest error observed during the parameter studies with the numerical error obtained for the discretization as predicted by the estimate, it appears that -- in the worst case -- the error is $2\dots3$ times larger than the actual minimum. However, taking into account the high sensitivity of the numerical error with respect to different values of the dimensionless Young's modulus, this deviation is still acceptable. In summary, it can be concluded that the error estimate from Section \ref{sec:conv} provides adequate guidance for a beneficial discretization for a broad range of Poisson's ratios.

Additionally, the results in Figure \ref{fig:num_leadingError} confirm another observation from the leading-order error analysis, which predicted that the minimal achievable error increases when moving closer to incompressible material behavior. Comparing the numerical errors for the cases with $\nu=0.1$ and $\nu=0.9$ shows that when choosing the best-possible combination of discretization parameters in each case, the error for $\nu=0.9$ is roughly one order of magnitude larger. As a result, finer meshes are required to achieve the same solution accuracy with larger Poisson's ratios.

\subsubsection{Stability boundaries in the material parameter space}

The next numerical study investigates the stability bounds of Section \ref{sec:stability}. To this end, several examples with material parameters close to the stability boundary are simulated. The examples are grouped into pairs, where for each pair one combination of discretization parameters lies in the stable domain and the other in the unstable regime. For the material parameters used in this investigation see Figure \ref{fig:stability}.

In order to ensure that the theoretically predicted onset of instability is actually triggered during the numerical experiments, a manufactured solution with a dense frequency spectrum is selected. Accordingly, the target solution is a 2D Gaussian hill located in the center of the problem domain $[0,L]\times[0,L]$ 
\begin{align}
    \label{eq:gaussian_x}
    \hat{u}_x/U = 10^{-3} \exp{\left[-\frac{(\tilde{x}-1/2)^2}{\sigma_1^2}\right]}\exp{\left[-\frac{(\tilde{y}-1/2)^2}{\sigma_2^2}\right]} - 10^{-3} \pi\sigma_1\sigma_2 \erf{(1/(2\sigma_1))}\erf{(1/(2\sigma_2))} \\
    \label{eq:gaussian_y}
    \hat{u}_y/U = 10^{-3} \exp{\left[-\frac{(\tilde{x}-1/2)^2}{\sigma_3^2}\right]}\exp{\left[-\frac{(\tilde{y}-1/2)^2}{\sigma_4^2}\right]} - 10^{-3} \pi\sigma_3\sigma_4 \erf{(1/(2\sigma_3))}\erf{(1/(2\sigma_4))}.
\end{align}
The constant offset ensures that the solution has zero mean. In the equation above, erf is the error function and the values of the standard deviations aligned with the coordinate directions $\sigma_i,i=1\dots4$ are listed in Table \ref{tab:std_gauss}.

\begin{table}[htbp]
    \centering
    \caption{Standard deviations $\sigma_i$ for the Gaussian hill manufactured solution}
    \begin{tabular}{r| *{4}{r}}
        $i$ & 1 & 2 & 3 & 4 \\ \hline &&&&\\[-1em]
        $\sigma_i$ & $120^{-1/2}$ & $123^{-1/2}$ & $118^{-1/2}$ & $125^{-1/2}$
    \end{tabular}
    \label{tab:std_gauss}
\end{table}

\begin{figure}[htbp]
    \centering
    \begin{subfigure}[b]{0.33\textwidth}
         \centering
         \includegraphics{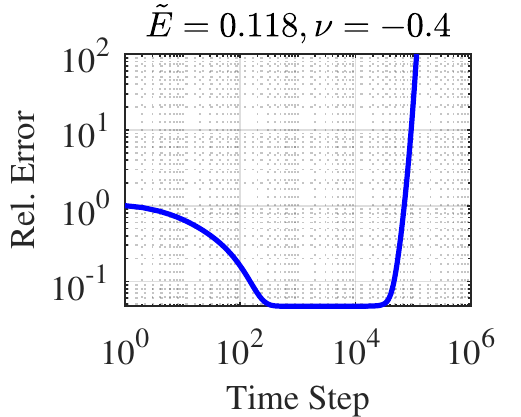}
         \caption{}
     \end{subfigure}
     \hfill
    \begin{subfigure}[b]{0.33\textwidth}
         \centering
         \includegraphics{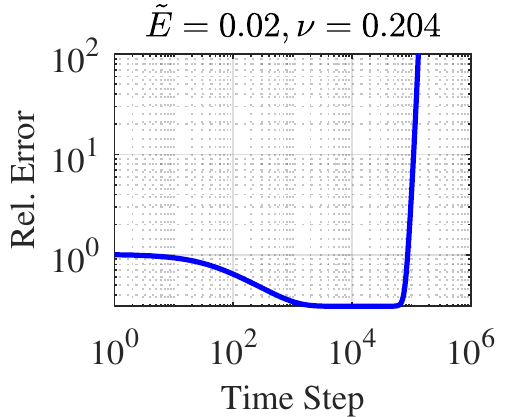}
         \caption{}
     \end{subfigure}
     \hfill
     \begin{subfigure}[b]{0.33\textwidth}
         \centering
         \includegraphics{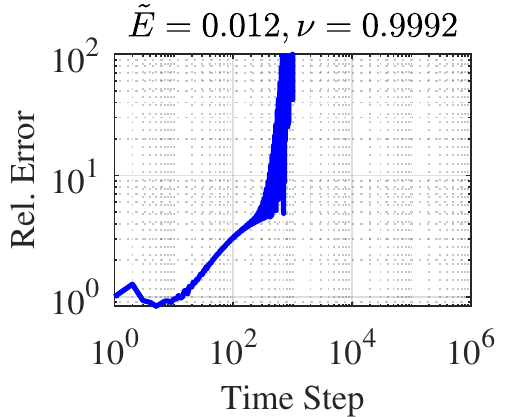}
         \caption{}
     \end{subfigure}
     
     \begin{subfigure}[b]{0.33\textwidth}
         \centering
         \includegraphics{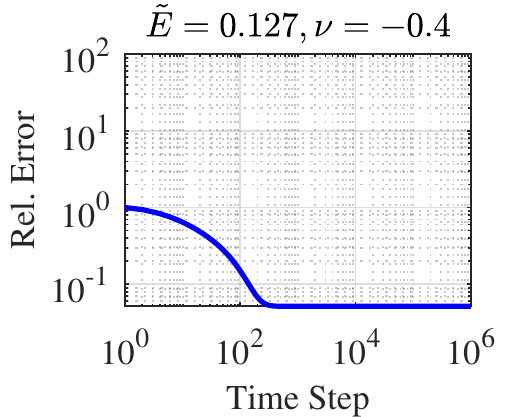}
         \caption{}
     \end{subfigure}
     \hfill
    \begin{subfigure}[b]{0.33\textwidth}
         \centering
         \includegraphics{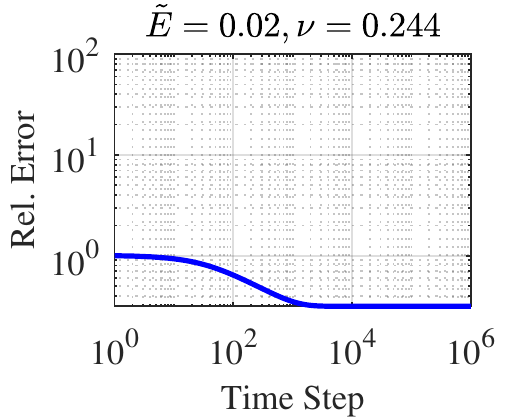}
         \caption{}
     \end{subfigure}
     \hfill
     \begin{subfigure}[b]{0.33\textwidth}
         \centering
         \includegraphics{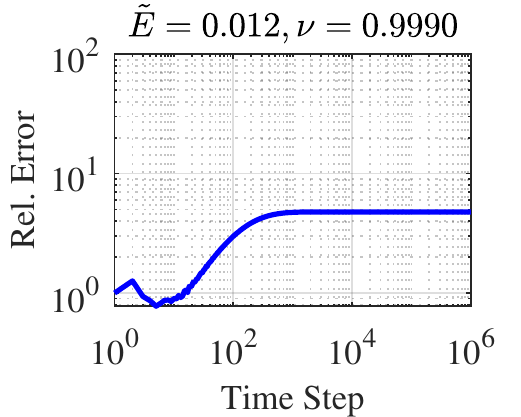}
         \caption{}
     \end{subfigure}
    \caption{Numerical investigation of the stability boundary showing the temporal evolution of the L2 norm of the displacement error. Top: Unstable discretization, bottom: Stable discretization. See Figure \ref{fig:stability} for the location of the dimensionless material parameters with respect to the stability bounds.}
    \label{fig:num_stability}
\end{figure}

The numerical examples in Figure \ref{fig:num_stability} demonstrate that the semi-analytically determined stability boundary is very accurately recovered by the given examples. The first two pairs of simulations (subfigures a), b), d) and e)) show that a high number of iterations is needed to actually observe the unstable behavior. In contrast, the last pair of simulations, which is carried out for nearly incompressible material behavior, displays a significantly faster onset of the instability once the theoretically predicted stability boundary is crossed. This differing behavior may result from the fact that the dimensionless material parameters are chosen in close vicinity to the stability boundary (see Figure \ref{fig:stability}) and that the onset of instability might not be as sharp as in theory or at exactly the predicted location because of numerical round-off or spectral filtering due to a finite mesh resolution \citep{Coreixas2020}.

Of the three considered pairs of material parameters, the stability region of the method poses an effective limitation only for $\nu=-0.4$. For this Poisson's ratio the leading-order error investigation suggests a good value of $\tilde{E}$ to be $0.082$, which lies well beyond the stability boundary that has been numerically verified in Figure \ref{fig:num_stability}. Note that this deterioration of the numerically achievable accuracy even increases for smaller values of the Poisson's ratio, because the advantageous choice of $\tilde{E}$ from the leading-order error analysis moves further into the unstable regime shown in Figure \ref{fig:stability}.

\subsubsection{Comparison with linear finite element results}

The final numerical study compares the results obtained by the improved lattice Boltzmann method (see Section \ref{sec:improved_scheme}) with simulations using standard bilinear quadrilateral  finite elements. For this comparison the more challenging solution of Eqs. \eqref{eq:gaussian_x} and \eqref{eq:gaussian_y} is considered. This manufactured solution results in a body load which is fairly concentrated in the interior of the domain, so that the problem can be considered to be approximately periodic. It is therefore acceptable to solve this problem with the lattice Boltzmann scheme using periodic boundary conditions (see Section \ref{sec:icbc}). For the finite element analysis it is more straightforward to employ Dirichlet-type boundary conditions that prescribe the analytical solution on the whole domain boundary.

\begin{figure}[htbp]
    \centering
    \begin{minipage}{0.5\textwidth}
    \centering
    \includegraphics{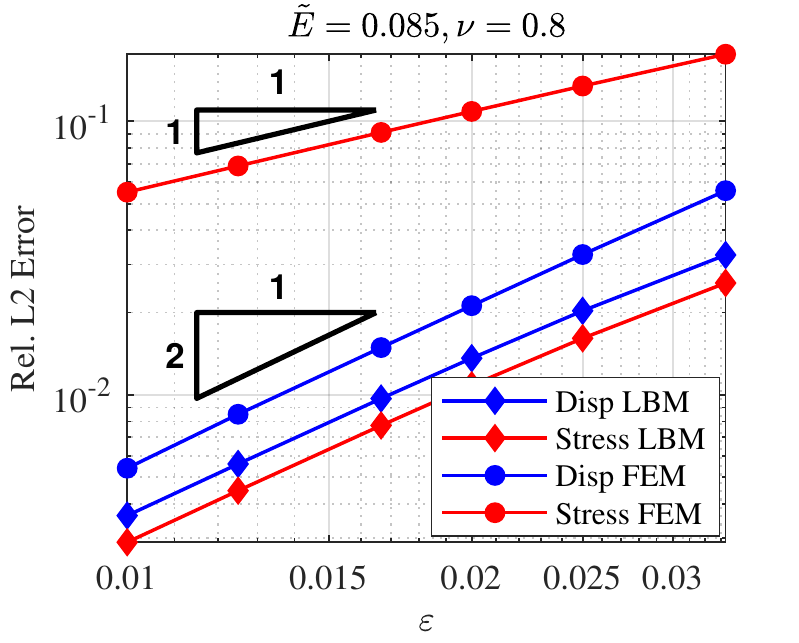}
    \end{minipage}
    \begin{minipage}{0.3\textwidth}
%    \begin{tabular}{l  l}
%    Norm & Conv. Order \\ \hline
%        L2 disp & 2.04 \\
%        Linf disp & 2.03 \\
%        L2 stress & 1.95 \\
%        Linf stress & 2.01
%    \end{tabular}
    \end{minipage}
    \caption{Comparison of the L2 errors using the improved lattice Boltzmann scheme and the finite element method}
    \label{fig:periodic_comparison}
\end{figure}

For both methods a convergence study is carried out using a structured mesh with element size $\Delta x = \varepsilon L$ in the case of the finite element analysis and a regular lattice with the same dimensionless grid spacing $\varepsilon$ for the lattice Boltzmann method. For the convergence study shown in Figure \ref{fig:periodic_comparison}, the number of (pseudo-)time steps to obtain the static equilibrium solution with the lattice Boltzmann method ranges between 5\,400 and 60\,000 steps and follows the diffusive scaling assumption, i.\,e. $\Delta t=\varepsilon^2T$ (see Section \ref{sec:asympt}).

The results of the convergence studies in Figure \ref{fig:periodic_comparison} reveal similar error behavior in the L2 norm for the displacement solution. Overall, slightly smaller errors in the displacement solution can be observed for the lattice Boltzmann method. However, as expected from the discussion in Section \ref{sec:stress}, the L2 errors in the Cauchy stress approximation decrease one order faster for the lattice Boltzmann scheme because of the higher consistency order, and are also significantly smaller. In summary, when comparing the accuracy on analogous discretizations, the lattice Boltzmann method performs similarly well in the displacement solution, but significantly better for the Cauchy stress.

\section{Conclusion}

We developed a novel lattice Boltzmann scheme to solve the quasi-static equations of linear elasticity. Because of the explicit time stepping inherent to the method, the target equation is extended by a damping term. Accordingly, the static equilibrium solution is obtained at the end of a transient phase that starts from a given initial condition.
The proposed method is devised using the generalized multiple relaxation time collision operator as a starting point. This allows to independently adjust the two material parameters of linear isotropic elasticity. In contrast to previous works, three important properties of the method are established:
\begin{enumerate}
    \item Only a single distribution function is required to determine the vector-valued solution field;
    \item Only a standard velocity set without rest population is necessary, i.\,e. the D2Q8 set in 2D;
    \item The scheme involves no finite difference approximations.
\end{enumerate}
As a result, the computational efficiency and algorithmic simplicity of the native lattice Boltzmann method are fully retained for this new application.

Using the asymptotic expansion technique, second-order consistency is demonstrated analytically and a leading-order error analysis provides effective guidance on how to adjust the numerical parameters of grid spacing, time step size and damping coefficient for improved accuracy. A stability analysis of the linear operator reveals a very broad range of Poisson's ratios that can be handled by the method. All analytical results are in excellent agreement with numerical results obtained with the method of manufactured solutions.

This initial contribution was limited only to the case of periodic problems in order to establish a thorough understanding of the method in the bulk. In a forthcoming contribution, the scheme will be furnished with suitable boundary formulations to be able to simulate physically relevant problems with Dirichlet- and Neumann-type boundaries.

With respect to established numerical methods for the solution of the linear elasticity equation, the lattice Boltzmann method proposed in this paper has drawbacks (mainly the need for time stepping to recover the quasi-static solution at steady state) but also advantages (second-order consistency also in the stresses, algorithmic simplicity and good scaling in parallelization). Based on the first results obtained in this paper, we believe that the method holds potential for further development. On the other hand, the present contribution may serve as an initial step and foundation for future investigations to solve more complex (nonlinear and/or dynamic) problems. Finally, in the context of multi-physics problems - such as for the simulation of additive manufacturing processes - it is anticipated that the handling of all physics (e.\,g. solid and fluid mechanics) within a single computational framework comes with significant advantages over the often used alternative of coupling separate codes based on different methods.

\section*{Declaration of competing interest}
The authors declare that they have no known competing financial interests or personal relationships that could have appeared to influence the work reported in this paper.

\section*{CRediT authorship contribution statement}
\textbf{Oliver Boolakee}: Conceptualization, Methodology, Software, Formal analysis, Investigation, Data curation, Writing – original draft, Writing – Review \& Editing, Visualization. \textbf{Martin Geier}: Methodology, Formal analysis, Investigation, Writing – review \& editing. \textbf{Laura De Lorenzis}: Conceptualization, Methodology, Writing – review \& editing, Supervision, Funding acquisition.

\section*{Acknowledgments}
Martin Geier acknowledges financial support by the German Research Foundation (DFG) project number 414265976-TRR 277.

\appendix

\section{Asymptotic expansion examples} \label{sec:expansion_examples}

In order to facilitate the understanding of the asymptotic expansion expression in Eq. \eqref{eq:expansion} a few examples with concrete values for $a,\:b$ and $r$ are provided without entering the recursion so that the results in Section \ref{sec:asympt} can be more easily reproduced. Note that the aliasing properties of the D2Q8 velocity set (e.\,g. $m_{30}=m_{10}$ or $m_{13}=m_{11}$) have already been applied in the expressions below.

\begin{itemize}
    \item $r=1$
\end{itemize}
\begin{align}
    a=1,b=0: \qquad \Omega_{10}^{(1)} = \partial_{\tilde{x}} \bar{m}_{20}^{(0)} + \partial_{\tilde{y}} \bar{m}_{11}^{(0)} \\
    a=1,b=1: \qquad \Omega_{11}^{(1)} = \partial_{\tilde{x}} \bar{m}_{21}^{(0)} + \partial_{\tilde{y}} \bar{m}_{12}^{(0)} \\
    a=2,b=0: \qquad \Omega_{20}^{(1)} = \partial_{\tilde{x}} \bar{m}_{10}^{(0)} + \partial_{\tilde{y}} \bar{m}_{21}^{(0)} \\
    a=0,b=2: \qquad \Omega_{02}^{(1)} = \partial_{\tilde{x}} \bar{m}_{12}^{(0)} + \partial_{\tilde{y}} \bar{m}_{01}^{(0)}
\end{align}

\begin{itemize}
    \item $r=2$
\end{itemize}
\begin{align}
    a=1,b=0: \qquad \Omega_{10}^{(2)} = \partial_{\tilde{t}} \bar{m}_{10}^{(0)} + \partial_{\tilde{x}} \bar{m}_{20}^{(1)} + \partial_{\tilde{y}} \bar{m}_{11}^{(1)} \\
    a=1,b=1: \qquad \Omega_{10}^{(2)} = \partial_{\tilde{t}} \bar{m}_{11}^{(0)} + \partial_{\tilde{x}} \bar{m}_{21}^{(1)} + \partial_{\tilde{y}} \bar{m}_{12}^{(1)} \\
    a=2,b=0: \qquad \Omega_{20}^{(2)} = \partial_{\tilde{t}} \bar{m}_{20}^{(0)} + \partial_{\tilde{x}} \bar{m}_{10}^{(1)} + \partial_{\tilde{y}} \bar{m}_{21}^{(1)} \\
    a=0,b=2: \qquad \Omega_{02}^{(2)} = \partial_{\tilde{t}} \bar{m}_{02}^{(0)} + \partial_{\tilde{x}} \bar{m}_{12}^{(1)} + \partial_{\tilde{y}} \bar{m}_{01}^{(1)}
\end{align}

\section{Constants of the leading-order error coefficients} \label{sec:explicit_constants}

For completeness the explicit expressions of all constants appearing in the body force term of the governing equation of the leading-order error coefficient are listed below (see Eq. \eqref{eq:leadingCoefficients} in Section \ref{sec:conv}). The relaxation times $\tau_\alpha$ of the second-order moments have been replaced with the dimensionless material parameters using Eqs. \eqref{eq:match11}--\eqref{eq:matchs}. In order to see how the time step size $\Delta t$ and the damping parameter $\kappa$ affect the dimensionless Young's modulus see Eq. \eqref{eq:scalingE}.

\begin{align}
    C_1 = &\frac{\nu^2-2\theta\nu+1}{(1-\theta^2)(1-\nu^2)^3}\tilde{E}^3 + \frac{\theta(\nu-\theta)\tau_{12}}{(1-\theta^2)(1-\nu^2)^2}\tilde{E}^2 + \frac{-1+\theta^2+12\theta(\theta-\nu)\tau_{12}\tau_{22}}{12(1-\theta^2)(1-\nu^2)}\tilde{E} \\
    C_2 = & \frac{(1-\nu)^3+4\theta(\nu^3-\nu^2+3\nu+1)+\theta^2(5\nu^3-15\nu^2-\nu-5)}{8(1-\theta^2)(1-\nu^2)^3}\tilde{E}^3 \nonumber \\ &+ \frac{2\theta\left((1-\nu)\nu+\theta^2(1+\nu)-\theta(1+\nu^2)\right)\tau_{12}+(1-\theta)(1-\nu)(1+2\theta-\nu-\theta^2(1+\nu))\tau_{21}}{-4\theta(1-\theta^2)(1-\nu^2)^2}\tilde{E}^2 \nonumber \\
    &+ \frac{(2\theta^2-3\theta+1)(1+3\theta(1-\nu)-5\nu)+24\theta(\theta-\nu)\tau_{12}\tau_{22}}{24(1-\theta^2)(1-\nu^2)}\tilde{E} \\
    C_3 = &\frac{(1-\nu)^2(3+\nu) + 4\theta(-\nu^3+3\nu^2+\nu+1)-\theta^2(5\nu^3-3\nu^2+7\nu+7)}{8\theta(1-\theta^2)(1-\nu^2)^3}\tilde{E}^3 \nonumber \\ &+\frac{-2(1-\theta)\theta(\nu-2\theta-1)(1-\nu)\tau_{12}-\left( (1+\theta)(1-\nu)^2+\theta^3(1+\nu)^2+\theta^2(\nu^2-2\nu-3) \right)\tau_{21}}{4\theta(1-\theta^2)(1-\nu^2)^2}\tilde{E}^2 \nonumber \\ &+\frac{(1-\theta)\left(1-\nu-\theta(1+\nu)-2\theta^2(2+\nu)\right)-4(1-\theta^2)(1-\nu)\tau_{12}\tau_{22}+4(\theta^2(1+\nu)+\nu-2\theta-1)\tau_{21}\tau_{22}}{8(1-\theta^2)(1-\nu^2)}\tilde{E} \nonumber \\ &+ \frac{\theta\tau_{22}}{4} \\
    C_4 = & \frac{16\theta\nu+(1-\nu)^2(3+\nu)-\theta^2(\nu^3+9\nu^2-5\nu+11)}{8\theta(1-\theta^2)(1-\nu^2)^3}\tilde{E}^3 \nonumber \\ &+\frac{(1-\theta^2)(1-\nu)(\theta\nu+\nu+\theta-1)\tau_{12}-2\theta\left(\theta^2(1+\nu)+\nu-2\theta-1\right)\tau_{21}}{4\theta(1-\theta^2)(1-\nu^2)^2}\tilde{E}^2 \nonumber \\ &+\frac{(1-\theta^2)\left( 2(2-\nu)-3\theta(5+\nu) \right)-12(1-\theta^2)(1-\nu)\tau_{12}\tau_{22}+12\left( \theta^2(1+\nu)+\nu-2\theta-1 \right)\tau_{21}\tau_{22}}{24(1-\theta^2)(1-\nu^2)}\tilde{E} \nonumber \\ &+\frac{\theta\tau_{22}}{4} \\
    C_5 = &\frac{1}{8\theta(1+\nu)^3}\tilde{E}^3-\frac{(1-\theta)\tau_{12}}{4\theta(1+\nu)^2}\tilde{E}^2 -\frac{3\theta-2}{24(1+\nu)}\tilde{E} \\
    D_1 = &\frac{\nu^2-2\theta\nu+1}{(1-\theta^2)(1-\nu^2)^2}\tilde{E}^2 + \frac{\theta(\nu-\theta)\tau_{12}}{(1-\theta^2)(1-\nu^2)^2}\tilde{E} - \frac{1}{4} \\
    D_2 = &\frac{(1-\nu)^2+8\theta\nu-\theta^2(5\nu^2-2\nu+5)}{4\theta(1-\theta^2)(1-\nu^2)^2}\tilde{E}^2 - \frac{\left( \theta^2(1+\nu)+\nu-2\theta-1\right)\tau_{21}}{2(1-\theta^2)(1-\nu^2)^2}\tilde{E}-\frac{\theta}{2} \\
    D_3 = &\frac{1}{4\theta(1+\nu)^2}\tilde{E}^2 + \frac{\tau_{12}}{2(1+\nu)}\tilde{E}-\frac{\theta}{4}
\end{align}

\section{Choice of relaxation rates for a fourth-order convergent scheme} \label{sec:rates_fourth_order}

The third convergence study in Section \ref{sec:num} shows an example with fourth-order convergence in the displacement field solution. In order to achieve this result, the exact solution needs to have vanishing mixed derivatives so that only two constants remain in the body force terms $r_x^{(4)}$ and $r_y^{(4)}$ of the governing equation of the leading-order error expansion. In order to fully remove this term, the following expressions must be used for the higher-order relaxation times.

\begin{align}
    \tau_{12}&=\tau_{21}=\frac{2-3\theta+12\theta\tau_{11}^2}{12(1-\theta)\tau_{11}} \\
    \tau_{22} &= \frac{\left( (1-\theta)\tau_d+(1+\theta)\tau_s \right) \left( \theta(2-3\theta)(\tau_s-\tau_d)-12\theta^2\tau_{11}^2(\tau_s-\tau_d)+2\tau_{11} \left( 1-\theta - 6(1-\theta)^2\tau_d^2-6(1-\theta^2)\tau_s^2 \right) \right)}{2\theta(2-3\theta+12\theta\tau_{11}^2)(\tau_s-\tau_d)}
\end{align}
Clearly, this has to be combined with the body load correction in Eqs. \eqref{eq:force_correction_x} and \eqref{eq:force_correction_y}.

%% References


\bibliographystyle{elsarticle-num}
\bibliography{bibliography}






% \end{thebibliography}

\end{document}